\documentclass[10pt,journal]{IEEEtran}

\usepackage{xspace,amssymb,epsfig,subfigure,syntonly}
\usepackage{epsfig,amsmath,color}
\usepackage{xspace,syntonly}
\usepackage{url}
\usepackage{algorithm,algorithmic}

\newcommand{\utwi}[1]{\mbox{\boldmath $#1$}}

\newcommand{\trace}{{\textrm{Tr}}}
\newcommand{\rank}{{\textrm{rank}}}

\newcommand{\cD}{{\cal D}}

\newcommand{\cL}{{\cal{L}}}
\newcommand{\cN}{{\cal N}}
\newcommand{\cP}{{\cal P}}

\newcommand{\cG}{{\cal G}}

\newcommand{\cR}{{\cal R}}

\newcommand{\cC}{{\cal C}}
\newcommand{\cE}{{\cal E}}

\newcommand{\cF}{{\cal F}}
\newcommand{\cB}{{\cal B}}
\newcommand{\cU}{{\cal U}}

\newcommand{\cH}{{\cal H}}

\newcommand{\ba}{{\bf a}}
\newcommand{\bb}{{\bf b}}

\newcommand{\be}{{\bf e}}

\newcommand{\bp}{{\bf p}}
\newcommand{\bq}{{\bf q}}

\newcommand{\bv}{{\bf v}}
\newcommand{\bi}{{\bf i}}

\newcommand{\bA}{{\bf A}}
\newcommand{\bB}{{\bf B}}

\newcommand{\bL}{{\bf L}}
\newcommand{\bM}{{\bf M}}

\newcommand{\bQ}{{\bf Q}}

\newcommand{\bI}{{\bf I}}
\newcommand{\bX}{{\bf X}}

\newcommand{\bW}{{\bf W}}
\newcommand{\bY}{{\bf Y}}
\newcommand{\bV}{{\bf V}}

\newcommand{\bPsi}{{\utwi{\Psi}}}

\newcommand{\bUpsilon}{{\utwi{\Upsilon}}}

\newcommand{\sfH}{\textsf{H}}
\newcommand{\sfT}{\textsf{T}}

\usepackage{epstopdf}
\usepackage[noadjust]{cite}

\begin{document}

\newtheorem{definition}{Definition}
\newtheorem{remark}{Remark}
\newtheorem{proposition}{Proposition}
\newtheorem{lemma}{Lemma}

\title{Decentralized Optimal Dispatch of  Photovoltaic Inverters in Residential Distribution Systems} 
\author{Emiliano Dall'Anese, \emph{Member}, \emph{IEEE}, Sairaj V. Dhople, \emph{Member}, \emph{IEEE}, Brian B. Johnson, \emph{Member}, \emph{IEEE},  \\ and Georgios B. Giannakis, \emph{Fellow}, \emph{IEEE}
\thanks{\protect\rule{0pt}{0.5cm}%

Submitted March 5, 2013; revised July 23, 2014; accepted September 4, 2014. 

This work was supported by NSF-CCF grants no. 1423316 and CyberSEES 1442686, the Institute of Renewable Energy and the Environment (IREE) grant no. RL-0010-13, University of Minnesota, and by the Laboratory Directed Research and Development (LDRD) Program at the National Renewable Energy Laboratory. 

E. Dall'Anese, S. Dhople, and G. Giannakis are with the  Dept. of ECE and Digital Technology Center, University of Minnesota, 200 Union Street SE, Minneapolis, MN, USA; e-mails: {\tt \{emiliano, sdhople, georgios\}@umn.edu}. B. Johnson is with the National Renewable Energy Laboratory, Golden, CO, USA; e-mail: {\tt brian.johnson@nrel.gov}  
}
}

\markboth{IEEE TRANSACTIONS ON ENERGY CONVERSION}
{Dall'Anese \MakeLowercase{\textit{et al.}}: }

\maketitle

\vspace{.5cm}

\begin{abstract}

  Decentralized methods for computing optimal real and reactive power
  setpoints for residential photovoltaic (PV) inverters are developed
  in this paper. It is known that conventional PV inverter
  controllers, which are designed to extract maximum power at unity
  power factor, cannot address secondary performance objectives such
  as voltage regulation and network loss minimization. Optimal power
  flow techniques can be utilized to select which inverters will
  provide ancillary services, and to compute their optimal real and
  reactive power setpoints according to well-defined performance
  criteria and economic objectives. Leveraging advances in
  sparsity-promoting regularization
  techniques and semidefinite relaxation, this paper shows how such problems can be solved with
  reduced computational burden and optimality guarantees. To enable
  large-scale implementation, a novel algorithmic framework is
  introduced --- based on the so-called alternating direction method
  of multipliers --- by which optimal power flow-type problems in this
  setting can be systematically decomposed into sub-problems that can
  be solved in a decentralized fashion by the utility and
  customer-owned PV systems with limited exchanges of
  information. Since the computational burden is shared among multiple
  devices and the requirement of all-to-all communication can be
  circumvented, the proposed optimization approach scales favorably to
  large distribution networks.

\end{abstract}


\section{Introduction}
\label{sec:Introduction}

\PARstart{T}{HE} proliferation of residential-scale photovoltaic (PV)
systems has highlighted unique challenges and concerns in the
operation and control of low-voltage distribution
networks. Secondary-level control of PV inverters can alleviate
extenuating circumstances such as overvoltages during periods when PV
generation exceeds the household demand, and voltage transients during
rapidly varying atmospheric conditions~\cite{Liu08}. Initiatives to
upgrade inverter controls and develop business models for ancillary
services are currently underway in order to
facilitate large-scale integration of renewables while ensuring
reliable operation of existing distribution feeders~\cite{Caramanis_PESTD14}.

  Examples of ancillary services include reactive power
  compensation, which has been recognized as a viable option to effect
  voltage regulation at the medium-voltage distribution
  level~\cite{Turitsyn11,Chertkov-ADMM13,Aliprantis13,Farivar12,Bolognani13}. The
  amount of reactive power injected or absorbed by inverters can be
  computed based on either local droop-type proportional
  laws~\cite{Turitsyn11,Aliprantis13}, or optimal power flow (OPF)
  strategies~\cite{Farivar12,Bolognani13}. Either way, voltage
  regulation with this approach comes at the expense of low power
  factors at the substation and high network currents, with the latter
  leading to high power losses in the
  network~\cite{Tonkoski11}. Alternative approaches require inverters
  to operate at unity power factor and to curtail part of the available
  active power~\cite{Tonkoski11,Samadi14}. For instance, heuristics
  based on droop-type laws are developed in~\cite{Tonkoski11} to
  compute the active power curtailed by each inverter in a residential
  system. Active power curtailment strategies are particularly
  effective in the low-voltage portion of distribution feeders, where
  the high resistance-to-inductance ratio of low-voltage overhead
  lines renders voltage magnitudes more sensitive to variations in the
  active power injections~\cite{kerstingbook}.

Recently, we proposed an optimal inverter dispatch (OID)
framework~\cite{OID} where the subset of critical PV-inverters that
most strongly impact network performance objectives are identified and
their real and reactive power setpoints are computed. This is
accomplished by formulating an OPF-type problem, which encapsulates
well-defined performance criteria as well as network and inverter
operational constraints. By leveraging advances in sparsity-promoting
regularizations and semidefinite relaxation (SDR) techniques~\cite{OID}, the problem
is then solved by a centralized computational device with reduced
computational burden. The proposed OID framework provides increased
flexibility over Volt/VAR
approaches~\cite{Turitsyn11,Aliprantis13,Farivar12,Bolognani13} and
active power curtailment methods~\cite{Tonkoski11,Samadi14} by: i) determining
in real-time those inverters that must participate in ancillary
services provisioning; and, ii) jointly optimizing both the real and
reactive power produced by the participating inverters (see, e.g.,
Figs.~\ref{Fig:OIDregions} (c)-(d) for an illustration of the the
inverters' operating regions under OID).

As proposed originally, the OID task can be carried out on a
centralized computational device which has to communicate with all
inverters. In this paper, the OID problem proposed
  in~\cite{OID} is strategically decomposed into sub-problems that can
  be solved in a decentralized fashion by the utility-owned energy
  managers and customer-owned PV systems, with limited exchanges of
  information. Hereafter, this suite of decentralized optimization
  algorithms is referred to as \emph{decentralized optimal inverter
    dispatch} (DOID). Building on the concept of leveraging both real
  and reactive power optimization~\cite{OID}, and decentralized
  solution approaches for OPF problems~\cite{Dallanese-TSG13}, two
  novel decentralized approaches are developed in this paper. In the
first setup, all customer-owned PV inverters can communicate with the
utility. The utility optimizes network performance (quantified in
terms of, e.g., power losses and voltage regulation) while individual
customers maximize their economic objectives (quantified in terms of,
e.g., the amount of active power they might have to curtail). This
setup provides flexibility to the customers to specify their
optimization objectives since the utility has no control on customer
preferences. In the spirit of the advanced metering infrastructure
(AMI) paradigm, utility and customer-owned EMUs exchange relevant
information~\cite{Samadi-SGComm10,GatsisTSG12} to \emph{agree} on the
optimal PV-inverter setpoints.  Once the
  decentralized algorithms have converged, the active and reactive
  setpoints are implemented by the inverter controllers. In the second DOID
approach, the distribution network is partitioned into
\emph{clusters}, each of which contains a set of customer-owned PV
inverters and a single cluster energy manager (CEM). A decentralized
algorithm is then formulated such that the operation of each cluster
is optimized and with a limited exchange of voltage-related messages,
the interconnected clusters consent on the system-wide voltage profile.
The decentralized OID frameworks are developed by leveraging the 
alternating direction method of multipliers (ADMM)~\cite{BeT89,BoydADMoM}.
 
  Related works include~\cite{Baldick99}, where
  augmented Lagrangian methods  (related to ADMM) were employed to decompose non-convex
  OPF problems for transmission systems into per-area instances,
  and~\cite{Nogales03,Hug09}, where standard Lagrangian approaches
  were utilized in conjunction with Newton methods.  ADMM was
  utilized in~\cite{ErsegheOPFADMM} to solve non-convex OPF renditions in
  a decentralized fashion, and in~\cite{Magnusson14}, where successive
  convex approximation methods were utilized to deal with nonconvex
  costs and constraints. In the distribution systems context,
  semidefinite relaxations of the OPF problem for
  balanced systems were developed in~\cite{Tse12}, and solved via
  node-to-node message passing by using dual (sub-)gradient ascent-based
  schemes. Similar message passing is involved in the ADMM-based
  decentralized algorithm proposed in~\cite{Chertkov-ADMM13} where a
  reactive power compensation problem based on approximate power flow
  models is solved. SDR of the OPF task in three-phase
  unbalanced systems was developed in~\cite{Dallanese-TSG13}; the
  resultant semidefinite program was solved in a distributed fashion by using 
  ADMM. 

\begin{figure}[b]
\begin{center}
\includegraphics{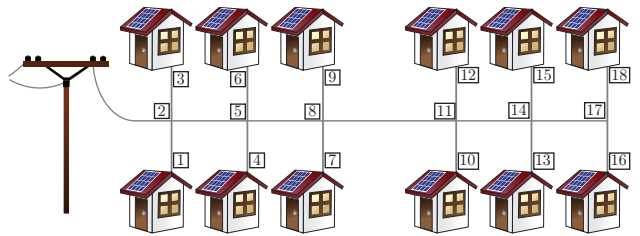}
\caption{Example of low-voltage residential network with high PV penetration, utilized in the test cases discussed in Section~\ref{sec:Simulations}. Node 0 corresponds to the secondary of the step-down transformer; set $\cU = \{2, 5, 8, 11, 14, 17\}$ collects nodes corresponding to distribution poles; and, homes $\mathrm{H}_1,\dots,\mathrm{H}_{12}$ are connected to nodes in the set $\cH = \{1,3,4,6,7,9,10,12,13,15,16,18\}$.}
\label{F_LV_network}
\vspace{-.5cm}
\end{center}
\end{figure}

The decentralized OID framework considerably broadens 
the setups of~\cite{Baldick99,Nogales03,Hug09,ErsegheOPFADMM,Magnusson14,Dallanese-TSG13,Tse12} by accommodating  different message passing strategies that are relevant in a variety of practical scenarios (e.g., customer-to-utility, customer-to-CEM and CEM-to-CEM communications). The proposed decentralized schemes offer improved optimality guarantees over~\cite{Baldick99,Nogales03,Hug09,ErsegheOPFADMM}, since it is grounded on an SDR technique; furthermore, ADMM enables superior convergence compared to~\cite{Tse12}. Finally, different from the distributed reactive compensation strategy of~\cite{Chertkov-ADMM13}, the proposed framework considers the utilization of an exact AC power flow model, as well as a joint computation of active and reactive power setpoint.   

For completeness, ADMM was utilized also in~\cite{Zhu_DSE,Kekatos_stateest13} for decentralized multi-area state estimation in transmission systems, and in~\cite{DallAneseTSTE14} to distribute over geographical areas the distribution system reconfiguration task.   

The remainder of the paper is organized as
follows. Section~\ref{sec:Formulation} briefly outlines
the centralized OID problem proposed
in~\cite{OID}. Sections~\ref{sec:Distributed}
and~\ref{sec:Distributed1} describe the two DOID problems discussed
above. Case studies to validate the approach are presented in
Section~\ref{sec:Simulations}. Finally, concluding remarks and
directions for future work are presented in
Section~\ref{sec:Conclusions}.

 \emph{Notation}:
  Upper-case (lower-case) boldface letters will be used for matrices
  (column vectors); $(\cdot)^\sfT$ for transposition; $(\cdot)^*$
  complex-conjugate; and, $(\cdot)^\sfH$ complex-conjugate
  transposition; $\Re\{\cdot\}$ and $\Im\{\cdot\}$ denote the real and
  imaginary parts of a complex number, respectively; $\mathrm{j} :=
  \sqrt{-1}$ the imaginary unit. $\trace(\cdot)$ the matrix trace;
  $\rank(\cdot)$ the matrix rank; $|\cdot|$ denotes the magnitude of a
  number or the cardinality of a set; $\|\bv\|_2 := \sqrt{\bv^\sfH
    \bv}$; $\|\bv\|_1 := \sum_i |[\bv]_i|$; and $\|\cdot\|_F$ stands
  for the Frobenius norm. Given a given matrix $\bX$, $[\bX]_{m,n}$ 
  denotes its $(m,n)$-th entry. 
  Finally, $\bI_N$ denotes the $N \times N$
  identity matrix; and, $\mathbf{0}_{M\times N}$, $\mathbf{1}_{M\times
    N}$ the $M \times N$ matrices with all zeroes and ones,
  respectively.

\section{Centralized optimal inverter dispatch} \label{sec:Formulation}

\subsection{Network and PV-inverter models}
\label{sec:inverter}

Consider a distribution system comprising $N+1$ nodes collected in the
set $\cN := \{0,1,\ldots,N\}$ (node $0$ denotes the secondary
of the step-down transformer), and lines represented by the set of
edges $\cE := \{(m,n)\} \subset \cN \times \cN$. For simplicity of exposition, a balanced system is
considered; however, both the centralized and decentralized frameworks
proposed subsequently can be extended to unbalanced systems
following the methods in~\cite{Dallanese-TSG13}. Subsets
$\cU, \cH \subset \cN$ collect nodes corresponding to utility poles
(with zero power injected or consumed), and those with installed
residential PV inverters, respectively (see Fig.~\ref{F_LV_network}).

\begin{figure}[b]
\begin{center}
\includegraphics[width=0.25\textwidth]{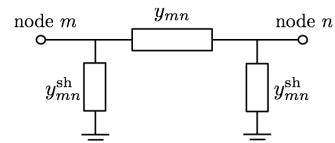}
\caption{$\pi$-equivalent circuits of a low-voltage distribution line $(m,n) \in \cE$.}
\label{F_line}
\vspace{-.5cm}
\end{center}
\end{figure}

Let $V_n\in \mathbb{C}$ and $I_n \in \mathbb{C}$ denote the phasors
for the line-to-ground voltage and the current injected at node $n \in
\cN$, respectively, and define $\bi := [I_1, \ldots, I_N]^\sfT \in
\mathbb{C}^{N}$ and $\bv := [V_1, \ldots, V_N]^\sfT \in
\mathbb{C}^{N}$. Using Ohm's and Kirchhoff's circuit laws, the linear
relationship $\bi = \bY \bv$ can be established, where the system
admittance matrix $\bY \in \mathbb{C}^{N+1 \times N+1}$ is formed
based on the system topology and the $\pi$-equivalent circuit of lines
$(m,n) \in \cE$, as illustrated in Fig.~\ref{F_line};
  see also~\cite[Chapter 6]{kerstingbook} for additional details on
  line modeling. Specifically, with $y_{mn}$ and
  $y^{\textrm{sh}}_{mn}$ denoting the series and shunt admittances of
  line $(m,n)$, the entries of $\bY$ are defined as
\begin{equation}
[\bY]_{m,n} := \left \{\begin{array}{ll} 
\sum_{j \in \cN_m} y^{\textrm{sh}}_{mj} + y_{mj} , & \textrm{if } m=n\\ 
- y_{mn} , &\textrm{if }(m,n)\in \cE\\
0, & \textrm{otherwise} \end{array} \right.
\nonumber
\end{equation}
where $\cN_{m} := \{j \in \cN: (m,j) \in \cE\}$ denotes the set of
nodes connected to the $m$-th one through a distribution line. 

\begin{figure}[t]
\begin{center}
\subfigure[]{\includegraphics{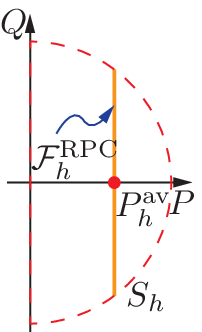}} 
\subfigure[]{\includegraphics{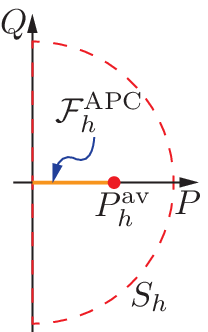}}
\subfigure[]{\includegraphics{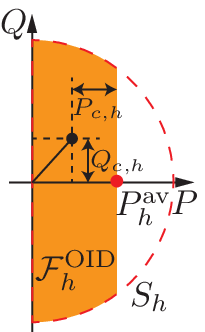}}
\subfigure[]{\includegraphics{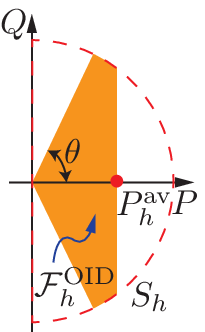}}
\end{center}
\caption{Feasible operating regions for the $h^\mathrm{th}$ inverter with apparent
  power rating $S_h$ under a) reactive power control (RPC), b) active power curtailment (APC), c) OID with joint control of real and reactive power,
and d) OID with a lower-bound on power factor.}
\label{Fig:OIDregions}
\end{figure}

A constant $PQ$ model~\cite{kerstingbook} is adopted for the load,
with  $P_{\ell,h}$ and $Q_{\ell,h}$ denoting the active and reactive
demands at node $h \in \cH$, respectively. For given solar irradiation
conditions, let $P_h^{\textrm{av}}$ denote the maximum
\emph{available active power} from the PV array at node $h \in
\cH$. The proposed framework calls for the joint control of
\emph{both} real and reactive power produced by the PV inverters. In
particular, the allowed operating regime on the complex-power plane
for the PV inverters is illustrated in Fig.~\ref{Fig:OIDregions} (d) and described by
\begin{equation}
\label{mg-PV} 
\cF^\mathrm{OID}_h := \left\{ P_{c,h}, Q_{c,h}:  
\begin{array}{l}
0  \leq P_{c,h}  \leq  P_{h}^{\textrm{av}}   \\
Q_{c,h}^2  \leq  S_{h}^2 - (P_{h}^{\textrm{av}}  - P_{c,h})^2 \\
| Q_{c,h} | \leq \tan \theta (P_{h}^{\textrm{av}}  - P_{c,h}) 
\end{array}
\hspace{-.2cm} \right\},
\end{equation}
where $P_{c,h}$ is the active power curtailed, and $Q_{c,h}$ is the
reactive power injected/absorbed by the inverter at node $h$. Notice
that if there is no limit to the power factor, then $\theta=\pi/2$,
and the operating region is given by Fig.~\ref{Fig:OIDregions}(c).

\subsection{Centralized optimization strategy}
\label{sec:Dispatch}

The centralized OID framework in~\cite{OID} invokes joint
optimization of active and reactive powers generated by the PV
inverters, and it offers the flexibility of selecting the subset of
critical PV inverters that should be dispatched in order to fulfill
optimization objectives and ensure electrical network constraints. To
this end, let $z_h$ be a binary variable indicating whether PV
inverter $h$ provides ancillary services or not and assume that at most $K < |\cH|$ PV
inverters are allowed to provide ancillary services. Selecting a (possibly time-varying)
subset of inverters promotes user fairness~\cite{Turitsyn11}, prolongs
inverter lifetime~\cite{Turitsyn11}, and captures possible
fixed-rate utility-customer pricing/rewarding strategies~\cite{Caramanis_PESTD14}. Let
$\bp_{c}$ and $\bq_{c}$ collect the active powers curtailed and the
reactive powers injected/absorbed by the inverters. With these
definitions, the OID problem is formulated as follows:
\begin{subequations} 
\label{Poidnonconvex}
\begin{align} 
\hspace{1.8cm} & \hspace{-1.7cm} \min_{\bv, \bi, \bp_{c}, \bq_{c},\{z_h\}} \,\,  C(\bV, \bp_{c})  \label{nonc-cost} \\
\mathrm{subject\,to} \,\, & \bi = \bY \bv, \, \{z_h\} \in \{0,1\}^{|\cH|}  \mathrm{~,and}  \nonumber  \\ 
V_h I_h^* & = (P_{h}^{\textrm{av}} - P_{c,h} - P_{\ell,h} ) + \mathrm{j}(Q_{c,h} - Q_{\ell,h}) \hspace{-.2cm} \label{nonc-balance} \\ 
 V_n I_n^* &  = 0    \hspace{3.2cm}  \forall \, n \in \cU \label{nonc-balancepoles}  \\
V^{\mathrm{min}} & \leq |V_n| \leq V^{\mathrm{max}}  \hspace{1.55cm}  \forall \, n \in \cN   \label{nonc-Vlimits} \\
(P_{c,h}, Q_{c,h}) & \in  \left\{
\begin{array}{l}
\{(0,0)\} , \mathrm{~if~} z_h = 0 \\
 \cF^\mathrm{OID}_h ,  \mathrm{~~~if~} z_h = 1
\end{array} 
\right.  \forall \, h \in \cH \label{nonc-oid} \\
\sum_{h \in \cH} z_h & \leq K, \,  \label{nonc-sn}
\end{align}
\end{subequations}
where constraint~\eqref{nonc-balance} is enforced at each node $h \in
\cH$; $C(\bV, \bp_{c})$ is a given cost function capturing both
network- and customer-oriented
objectives~\cite{OID,Caramanis_PESTD14};
and,~\eqref{nonc-oid}-\eqref{nonc-sn} jointly indicate which inverters
have to operate either under OID (i.e., $ (P_{c,h}, Q_{c,h})  \in
\cF^\mathrm{OID}_h$), or, in the business-as-usual mode (i.e., $ (P_{c,h},
Q_{c,h}) = (0,0)$). An alternative problem formulation can be 
obtained by removing constraint~\eqref{nonc-sn}, and adopting the cost $C(\bV, \bp_{c}) + \lambda_z \sum_{h \in \cH} z_h$
in~\eqref{nonc-cost}, with $\lambda_z \geq 0$ a weighting coefficient utilized  
to trade off achievable cost $C(\bV, \bp_{c})$  for the number of controlled inverters. When $\lambda_z$ represents a fixed reward for customers providing ancillary services~\cite{Caramanis_PESTD14} and $C(\bV, \bp_{c})$ models costs associated with active power losses and active power set points, OID~\eqref{Poidnonconvex} returns the inverter setpoints that minimize the economic cost incurred by feeder operation. 

As with various OPF-type problem formulations, the power balance
and lower bound on the voltage magnitude
constraints~\eqref{nonc-balance},~\eqref{nonc-balancepoles}
and~\eqref{nonc-Vlimits}, respectively, render the
OID problem~\emph{nonconvex}, and thus challenging to solve optimally
and efficiently. Unique to the OID formulation are the binary
optimization variables $\{z_h\}$; finding the optimal (sub)set of inverters to dispatch 
involves the solution of combinatorially many subproblems. 
Nevertheless, a computationally-affordable
\emph{convex} reformulation was developed in~\cite{OID}, by leveraging
 sparsity-promoting regularization~\cite{Wiesel11} and
semidefinite relaxation (SDR) techniques~\cite{LavaeiLow, Zhu_DSE,
  Dallanese-TSG13 } as briefly described next.

In order to bypass binary selection variables, key is to notice that that if 
inverter $h$ is \emph{not} selected for ancillary services, then one
clearly has that $P_{c,h} = Q_{c,h} = 0$ [cf.~\eqref{nonc-oid}]. Thus,
for $K < |\cH|$, one has that the $2 |\cH|\times 1$ real-valued vector
$[P_{c,1}, Q_{c,1}, \ldots, P_{c,|\cH|}, Q_{c,|\cH|} ]^\sfT$ is
\emph{group sparse}~\cite{Wiesel11}; meaning that, either the $2
\times 1$ sub-vectors $[P_{c,h}, Q_{c,h}]^\sfT$ equal  $\mathbf{0}$
or not~\cite{OID}. This group-sparsity attribute enables discarding the binary
variables and to effect PV inverter selection by regularizing the
cost in~\eqref{Poidnonconvex} with the following function:
\begin{align}
G(\bp_{c},\bq_{c})  :=  \lambda \sum_{h \in \cH} \,  \|[P_{c,h}, Q_{c,h}]^\sfT\|_2,  \label{Glasso_powers}
\end{align}
where $\lambda \geq 0$ is a tuning parameter. Specifically, the number
of inverters operating under OID decreases as $\lambda$ is
increased~\cite{Wiesel11}.

Key to developing a relaxation of the OID task is to
express powers and voltage magnitudes as linear functions of the
outer-product Hermitian matrix $\bV := \bv \bv^\sfH$, and to
reformulate the OID problem with cost and constraints that are linear
in $\bV$, as well as the constraints $\bV \succeq \mathbf{0}$ and
$\rank(\bV) = 1$~\cite{LavaeiLow,Zhu_DSE, Dallanese-TSG13}. The
resultant problem is still nonconvex because of the constraint
$\rank(\bV) = 1$; however in the spirit of SDR, this constraint can be
dropped. 

To this end, define the matrix $\mathbf{Y}_n := \be_n \be_n^\sfT
\mathbf{Y}$ per node $n$, where $\{\mathbf{e}_{n}\}_{n \in \cN}$
denotes the canonical basis of $\mathbb{R}^{|\cN|}$. Based on
$\mathbf{Y}_n$, define also the Hermitian matrices $\bA_{n} :=
\frac{1}{2} (\bY_n + \bY_n^\sfH)$, $\bB_{n} := \frac{j}{2} (\bY_n -
\bY_n^\sfH) $, and $\bM_{n} := \be_n \be_n^\sfT$. Using these
matrices, along with~\eqref{Glasso_powers} the relaxed convex OID
problem can be formulated as
\begin{subequations} 
\label{Poid}
\begin{align}
& \hspace{-1.9cm} \min_{\bV, \bp_{c}, \bq_{c}}   C(\bV, \bp_{c}) + G(\bp_{c},\bq_{c}) \label{Pm-cost}   \\
\,\, \mathrm{s.\,to} \,\,  \bV  \succeq \mathbf{0}, &  \mathrm{~and}  \nonumber  \\ 
\trace(\bA_h \bV)  & = - P_{\ell,h} + P_{h}^{\textrm{av}} - P_{c,h} \hspace{.65cm} \forall \, h \in \cH \label{Pm-balanceP}  \\
\trace(\bB_h \bV)  & = - Q_{\ell,h} + Q_{c,h}   \hspace{1.5cm} \forall \, h \in \cH \label{Pm-balanceQ} \\
\trace(\bA_n \bV)  &  = 0,~\trace(\bB_n \bV) = 0 \hspace{1.05cm} \forall \, n \in \cU  \label{Pm-balancePoles}  \\
 V_{\mathrm{min}}^2 & \leq \trace(\bM_n \bV)  \leq V_{\mathrm{max}}^2  \hspace{.85cm}  \forall \, n \in \cN   \label{Pm-Vlimits} \\
& \hspace{-.8cm} (P_{c,h}, Q_{c,h}) \in \cF^\mathrm{OID}_h \hspace{1.8cm} \forall \, h \in \cH. \label{Pm-oid}  
\end{align}
\end{subequations}
If the optimal solution of the relaxed problem~\eqref{Poid} has
rank 1, then the resultant voltages, currents, and power flows are globally optimal for
given inverter setpoints~\cite{LavaeiLow}.  Sufficient conditions for
SDR to be successful in OPF-type problems are available for
networks that are radial and balanced in~\cite{Lavaei_tree,Tse12},
whereas the virtues of SDR for unbalanced medium- and low-voltage
distribution systems have been demonstrated in~\cite{Dallanese-TSG13}.
As for the inverter setpoints $\{(P_{h}^{\textrm{av}} - P_{c,h}, Q_{c,h} )\}$, those
obtained from~\eqref{Poid} may be slightly sub-optimal compared to the setpoints that would have been obtained by 
solving the optimization problem~\eqref{Poidnonconvex}. This is mainly due to the so-called 
``shrinkage effect'' introduced by the regularizer~\eqref{Glasso_powers}~\cite{Wiesel11}.
Unfortunately, a numerical assessment of the optimality gap is impractical, since finding the globally optimal solution 
of problem~\eqref{Poidnonconvex} under all setups is  computationally infeasible.

To solve the OID problem, all customers' loads and available powers
$\{P_{h}^{\textrm{av}}\}$ must be gathered at a central processing
unit (managed by the utility company), which subsequently dispatches the PV
inverter setpoints. Next, decentralized implementations of
the OID framework are presented so that the OID problem can be solved in a decentralized fashion with limited exchange of information. From a computational perspective,
decentralized schemes ensure scalability of problem complexity with
respect to the system size.

\section{DOID: utility-customer message passing} \label{sec:Distributed}

Consider decoupling the cost $ C(\bV, \bp_{c})$ in~\eqref{Pm-cost} as
$ C(\bV, \bp_{c}) = C_{\textrm{utility}}(\bV, \bp_{c}) + \sum_h
R_{h}(P_{c,h})$, where $C_{\textrm{utility}}(\bV,\bp_{c})$ captures
utility-oriented optimization objectives, which may include e.g.,
power losses in the network and voltage
deviations~\cite{OID,Farivar12,Bolognani13}; and, $R_{h}(P_{c,h})$ is
a convex function modeling the cost incurred by (or the reward
associated with) customer $h$ when the PV inverter is required to
curtail power. Without loss of generality, a quadratic function
$R_{h}(P_{c,h}) := a_h P_{c,h}^2 + b_h P_{c,h}$ is adopted here, where
the choice of the coefficients is based on specific utility-customer
prearrangements~\cite{Caramanis_PESTD14} or customer preferences~\cite{OID}.

Suppose that customer $h$ transmits to the utility company the net active power $\bar{P}_{h}
:= - P_{\ell,h} + P^{\textrm{av}}_{h}$ and the reactive load
$Q_{\ell,h}$; subsequently, customer and utility will
agree on the PV-inverter setpoint, based on the optimization
objectives described by $C_{\textrm{utility}}$ and $\{R_{h}\}$. To
this end, let $\bar{P}_{c,h}$ and $ \bar{Q}_{c,h}$ represent \emph{copies} of
$P_{c,h}$, and $Q_{c,h}$, respectively, at the utility. The
corresponding $|\cH| \times 1$ vectors that collect the copies of the
inverter setpoints are denoted by $\bar{\bp}_{c}$ and $\bar{\bq}_{c}$,
respectively. Then, using the additional optimization variables
$\bar{\bp}_{c}, \bar{\bq}_{c}$, the relaxed OID problem~\eqref{Poid}
can be equivalently reformulated as:
\begin{subequations} 
\label{Pmconsensus1}
\begin{align}
& \hspace{-2.4cm} \min_{\substack{\bV, \bp_{c}, \bq_{c} \\ \bar{\bp}_{c},\bar{\bq}_{c}}} \,\,   \bar{C}(\bV, \bar{\bp}_{c},\bar{\bq}_{c}) + \sum_{h \in \cH} R_{h}(P_{c,h})  \label{Pmc-cost} \\
\,\, \mathrm{s.\,to} \,\,  \bV  \succeq \mathbf{0}, &  \mathrm{~and}  \nonumber  \\ 
\trace(\bA_h \bV)  & = \bar{P}_{h} - \bar{P}_{c,h} \hspace{2.05cm} \forall \, h \in \cH \label{Pmc-balanceP}  \\
\trace(\bB_h \bV)  & = - Q_{\ell,h} + \bar{Q}_{c,h}   \hspace{1.5cm} \forall \, h \in \cH \label{Pmc-balanceQ} \\
\trace(\bA_n \bV)  &  = 0,~\trace(\bB_n \bV) = 0 \hspace{1.05cm} \forall \, n \in \cU  \label{Pmc-balancePoles}  \\
 V_{\mathrm{min}}^2 & \leq \trace(\bM_h \bV)  \leq V_{\mathrm{max}}^2  \hspace{.85cm}  \forall \, n \in \cN   \label{Pmc-Vlimits} \\
& \hspace{-.8cm} (P_{c,h}, Q_{c,h}) \in \cF^\mathrm{OID}_h  \hspace{1.8cm} \forall \, h \in \cH \label{Pmc-oid}  \\
\bar{P}_{c,h} & = P_{c,h}, \,\ \bar{Q}_{c,h}  =  Q_{c,h}  \hspace{.8cm} \forall \, h \in \cH  \label{Pmc-consensusP} 
\end{align} 
\end{subequations}
where constraints~\eqref{Pmc-consensusP} ensure that utility and
customer \emph{agree} upon the inverters' setpoints, and
$\bar{C}(\bV, \bar{\bp}_{c},\bar{\bq}_{c}) :=  C_{\textrm{utility}}(\bV, \bar{\bp}_{c}) + G(\bar{\bp}_{c},\bar{\bq}_{c})$  
is the regularized cost function to be minimized at the utility.

The consensus constraints~\eqref{Pmc-consensusP} render problems~\eqref{Poid}
and~\eqref{Pmconsensus1} equivalent; however, the same constraints
impede problem decomposability, and thus modern optimization
techniques such as distributed (sub-)gradient
methods~\cite{Samadi-SGComm10,GatsisTSG12} and
ADMM~\cite[Sec.~3.4]{BeT89} cannot be directly applied to
solve~\eqref{Pmconsensus1} in a decentralized fashion. To enable problem
decomposability, consider introducing the auxiliary variables $x_h,
y_h$ per inverter $h$. Using these auxiliary
variables,~\eqref{Pmconsensus1} can be reformulated as
\begin{subequations} 
\label{Pmconsensus2}
\begin{align}
& \hspace{-2.4cm} \min_{\substack{\bV, \bp_{c}, \bq_{c} \\  \bar{\bp}_{c},\bar{\bq}_{c}, \{x_h, y_h\}}}   \bar{C}(\bV,  \bar{\bp}_{c},\bar{\bq}_{c}) + \sum_{h \in \cH} R_{h}(P_{c,h})  \label{mg-cost2} \\
\,\, \mathrm{s.\,to} \,\,  \bV  \succeq \mathbf{0}, &~\eqref{Pmc-balanceP}-\eqref{Pmc-oid},\mathrm{~and}  \nonumber  \\ 
\bar{P}_{c,h} & = x_h, \quad x_h = P_{c,h} \hspace{.5cm} \forall \, h \in \cH  \label{Pmc-consP2} \\
\bar{Q}_{c,h} & =  y_h, \quad y_h =  Q_{c,h}  \hspace{.5cm} \forall \, h \in \cH. \label{Pmc-consQ2} 
\end{align}
\end{subequations}
Problem~\eqref{Pmconsensus2} is equivalent to~\eqref{Poid}
and~\eqref{Pmconsensus1}; however, compared
to~\eqref{Poid}-\eqref{Pmconsensus1}, it is amenable to a decentralized
solution via ADMM~\cite[Sec.~3.4]{BeT89} as described in the remainder
of this section. ADMM is preferred over distributed (sub-)gradient
schemes because of its significantly faster
convergence~\cite{Chertkov-ADMM13} and resilience to communication
errors~\cite{ErsegheADMM}.

Per inverter $h$, let $\bar{\gamma}_h, \gamma_h$ denote the
multipliers associated with the two constraints in~\eqref{Pmc-consP2},
and $\bar{\mu}_h, \mu_h$ the ones associated
with~\eqref{Pmc-consQ2}. Next, consider the partial
quadratically-augmented Lagrangian of~\eqref{Pmconsensus2}, defined as
follows:
\begin{align}
& \cL(\bar{\cP}, \{\cP_h\}, \cP_{xy}, \cD )  := \bar{C}(\bV,  \bar{\bp}_{c},\bar{\bq}_{c}) + \sum_{h \in \cH} \Big[ R_{h}(P_{c,h}) \nonumber  \\ 
&  + \bar{\gamma}_h (\bar{P}_{c,h} - x_h) + \gamma_h(x_h - P_{c,h} ) + \bar{\mu}_h (\bar{Q}_{c,h} -  y_h) \nonumber  \\ 
& + \mu_h(y_h -  Q_{c,h}) +  (\kappa/2)(\bar{P}_{c,h} - x_h)^2  + (\kappa/2)(x_h - P_{c,h} )^2 \nonumber  \\ 
& + (\kappa/2)(\bar{Q}_{c,h} -  y_h)^2 + (\kappa/2)(y_h -  Q_{c,h})^2 \Big] \, ,
\label{eq:Lagrangian}
\end{align}
where $\bar{\cP} := \{\bV, \bar{\bp}_{c},\bar{\bq}_{c}\}$ collects the
optimization variables of the utility; $\cP_h := \{P_{c,h}, Q_{c,h}\}$
are the decision variables for customer $h$; $\cP_{xy} := \{x_h, y_h,
\forall h \in \cH \}$ is the set of auxiliary variables; $\cD :=
\{\bar{\gamma}_h, \gamma_h, \bar{\mu}_h, \mu_h, \forall h \in \cH \}$
collects the dual variables; and $\kappa > 0$ is a given
constant. Based on~\eqref{eq:Lagrangian}, ADMM amounts to iteratively
performing the steps \textbf{[S1]}-\textbf{[S3]} described next, where $i$ denotes the
iteration index:

\noindent \textbf{[S1]} Update variables $\bar{\cP}$ as follows:  
\begin{align}
\hspace{-.3cm} \bar{\cP}[i+1] := \arg &  \min_{\bV, \{\bar{P}_{c,h}, \bar{Q}_{s,h}\}} \cL(\bar{\cP}, \{\cP_h[i]\}, \cP_{xy}[i], \cD[i])  \label{ADMM1a} \\
& \,\, \mathrm{s.\,to} \,\,  \bV  \succeq \mathbf{0},\textrm{~and~} \eqref{Pmc-balanceP}-\eqref{Pmc-Vlimits} . \nonumber 
\end{align}
Furthermore, per inverter $h$, update $P_{c,h}, Q_{c,h}$ as follows: 
\begin{align}
\cP_h[i+1] := \arg & \hspace{-.3cm} \min_{P_{c,h}, Q_{c,h}} \cL(\bar{\cP}[i], P_{c,h}, Q_{c,h} , \cP_{xy}[i], \cD[i])  \label{ADMM1b} \\
& \,\, \mathrm{s.\,to} \,\, (P_{c,h}, Q_{c,h}) \in \cF^\mathrm{OID}_h  \nonumber 
\end{align}

\noindent \textbf{[S2]} Update auxiliary variables $\cP_{xy}$: 
\begin{align}
& \hspace{-.3cm} \cP_{xy}[i+1] := \nonumber \\
& \hspace{.3cm} \arg \min_{\{x_h, y_h\}} \cL(\bar{\cP}[i+1], \{\cP_h[i+1]\} , \{x_h, y_h\},  \cD[i])  \label{ADMM2} 
\end{align}

\noindent \textbf{[S3]} Dual update:  
\begin{subequations} 
\label{ADMM3}
\begin{align}
\bar{\gamma}_h[i+1] & =  \bar{\gamma}_h[i] + \kappa(\bar{P}_{c,h}[i+1] - x_h[i+1]) \label{ADMM3_bargam} \\
\gamma_h[i+1] & =  \gamma_h[i] + \kappa(x_h[i+1] - P_{c,h}[i+1]) \label{ADMM3_gam} \\
\bar{\mu}_h[i+1] & =  \bar{\mu}_h[i] + \kappa(\bar{Q}_{c,h}[i+1] -  y_h[i+1]) \label{ADMM3_barmu} \\
\mu_h[i+1] & =  \mu_h[i] + \kappa(y_h[i+1] - Q_{c,h}[i+1] ) . \label{ADMM3_mu}
\end{align}
\end{subequations}

\noindent In \textbf{[S1]}, the primal variables $\bar{\cP}, \{\cP_h\}$ are
obtained by minimizing~\eqref{eq:Lagrangian}, where the auxiliary
variables $\cP_{xy}$ and the multipliers $\cD$ are kept fixed to their
current iteration values. Likewise, the auxiliary variables are
updated in \textbf{[S2]} by fixing $\bar{\cP}, \{\cP_h\}$ to their up-to-date
values. Finally, the dual variables are updated in \textbf{[S3]} via dual
gradient ascent.

It can be noticed that step \textbf{[S2]} favorably decouples into $2|\cH|$
scalar and unconstrained quadratic programs, with $x_h[i+1]$ and
$y_h[i+1]$ solvable in closed-form. Using this feature, the
following lemma can be readily proved.

\begin{lemma}
\label{lemma:dualupdates}
Suppose that the multipliers are initialized as $\bar{\gamma}_h[0] = \gamma_h[0] = \bar{\mu}_h[0] = \mu_h[0] = 0$. Then, for all iterations $i>0$, it holds that: 

\noindent \emph{i)} $\bar{\gamma}_h[i] = \gamma_h[i]$; 

\noindent \emph{ii)} $\bar{\mu}_h[i] = \mu_h[i]$.  

\noindent \emph{iii)} $x_h[i] = \frac{1}{2}\bar{P}_{c,h}[i] + \frac{1}{2} P_{c,h}[i]$; and,

\noindent \emph{iv)} $y_h[i] = \frac{1}{2}\bar{Q}_{c,h}[i] + \frac{1}{2} Q_{c,h}[i]$. 

\end{lemma}

Using Lemma~\ref{lemma:dualupdates}, the conventional ADMM steps \textbf{[S1]}--\textbf{[S3]} can be simplified as follows.

\noindent \textbf{[S1$^\prime$]} At the utility side, variables
$\bar{\cP}$ are updated by solving the following convex optimization
problem:
\begin{subequations} 
\label{step1_utility} 
\begin{align}
\bar{\cP}[i+1] := \arg  & \min_{\bV, \{\bar{P}_{c,h}, \bar{Q}_{s,h}\}} \bar{C}(\bV, \bar{\bp}_{c},\bar{\bq}_{c})  \nonumber \\ 
&  \hspace{1.5cm} + F(\bar{\bp}_{c},\bar{\bq}_{c}, \{\cP_h[i]\}) \label{step1_utility_2} \\
& \,\, \mathrm{s.\,to} \,\,  \bV  \succeq \mathbf{0},\textrm{~and~} \eqref{Pmc-balanceP}-\eqref{Pmc-Vlimits} \nonumber 
\end{align}
where function $F(\bar{\bp}_{c},\bar{\bq}_{c}, \{\cP_h[i]\})$ is
defined as
\begin{align}
& F(\bar{\bp}_{c},\bar{\bq}_{c}, \{\cP_h[i]\}) := \sum_{h \in \cH} \Big[ \frac{\kappa}{2} (\bar{P}_{c,h}^2 + \bar{Q}_{c,h}^2) \nonumber \\
& \hspace{1.7cm} + \bar{P}_{c,h}\left(\gamma_h[i] - \frac{\kappa}{2} \bar{P}_{c,h}[i] - \frac{\kappa}{2} P_{c,h}[i] \right) \nonumber \\
& \hspace{1.7cm} + \bar{Q}_{c,h}\left(\mu_h[i] - \frac{\kappa}{2} \bar{Q}_{c,h}[i] - \frac{\kappa}{2} Q_{c,h}[i] \right) \Big] . \label{step1_utility_3}
\end{align}
\end{subequations}
At the customer side, the PV-inverter setpoints are updated by solving
the following constrained quadratic program:
\begin{align}
\cP_h[i+1] := \arg & \min_{P_{c,h}, Q_{c,h}} \Big[ R_h(P_{c,h}) + \frac{\kappa}{2}\left(P_{c,h}^2  + Q_{c,h}^2\right)  \nonumber \\
& \hspace{-.5cm} - P_{c,h}\left(\gamma_h[i] + \frac{\kappa}{2} \bar{P}_{c,h}[i] + \frac{\kappa}{2} P_{c,h}[i] \right) \nonumber \\
& \hspace{-.5cm} - Q_{c,h}\left(\mu_h[i] + \frac{\kappa}{2} \bar{Q}_{c,h}[i] + \frac{\kappa}{2} Q_{c,h}[i] \right) \Big]\label{step1_customer} \\
& \,\, \mathrm{s.\,to} \,\, (P_{c,h}, Q_{c,h}) \in \cF^\mathrm{OID}_h .  \nonumber 
\end{align}

\noindent \textbf{[S2$^\prime$]} At the utility and customer sides,
the dual variables are updated as:
\begin{subequations} 
\label{step2_customer_and_utility}
\label{step2}
\begin{align}
\gamma_h[i+1] & =  \gamma_h[i] + \frac{\kappa}{2} (\bar{P}_{c,h}[i+1] -  P_{c,h}[i+1])  \\
\mu_h[i+1] & =  \mu_h[i] + \frac{\kappa}{2} (\bar{Q}_{c,h}[i+1] -  Q_{c,h}[i+1]).
\end{align}
\end{subequations}

\begin{algorithm}[t]
\label{alg:OID2}
\caption{DOID: Utility-customer message passing} \small{
\begin{algorithmic}

\STATE Set $\gamma_h[0]  = \mu_h[0] = 0$ for all $h \in \cH$.

\FOR {$i = 1,2,\ldots$ (repeat until convergence)} 

\STATE 1. \textbf{[Utility]}: update $\bV[i+1]$ and $\{\bar{P}_{c,h}[i+1], \bar{Q}_{c,h}[i+1]\}$ via~\eqref{step1_utility}. \\

\hspace{.2cm}  \textbf{[Customer-$h$]}: update $\bar{P}_{c,h}[i+1], \bar{Q}_{c,h}[i+1]$ via~\eqref{step1_customer}. \\

\STATE 2.  \textbf{[Utility]}: send  $\bar{P}_{c,h}[i+1], \bar{Q}_{c,h}[i+1]$ to $h$; 

\hspace{1.5cm} repeat for all $h \in \cH$.

\hspace{.2cm}  \textbf{[Customer-$h$]}: receive $\bar{P}_{c,h}[i+1], \bar{Q}_{c,h}[i+1]$ from utility; 

\hspace{1.65cm} send $P_{c,h}[i+1], Q_{c,h}[i+1]$ to utility;

\hspace{1.65cm} repeat for all $h \in \cH$.

\hspace{.2cm}  \textbf{[Utility]}: receive $P_{c,h}[i+1], Q_{c,h}[i+1]$ from $h$;  

\hspace{1.5cm} repeat for all $h \in \cH$.

\STATE 3. \textbf{[Utility]}: update $\{\gamma_h[i+1] , \mu_h[i+1] \}_{h \in \cH}$ via~\eqref{step2}.

\hspace{.2cm}  \textbf{[Customer-$h$]}: update dual variables $\gamma_h[i+1] , \mu_h[i+1]$ via~\eqref{step2};

\hspace{1.65cm} repeat for all $h \in \cH$.

\ENDFOR

Implement setpoints in the PV inverters. 

\end{algorithmic}}
\end{algorithm}

\begin{figure}[t]
\begin{center}
\includegraphics{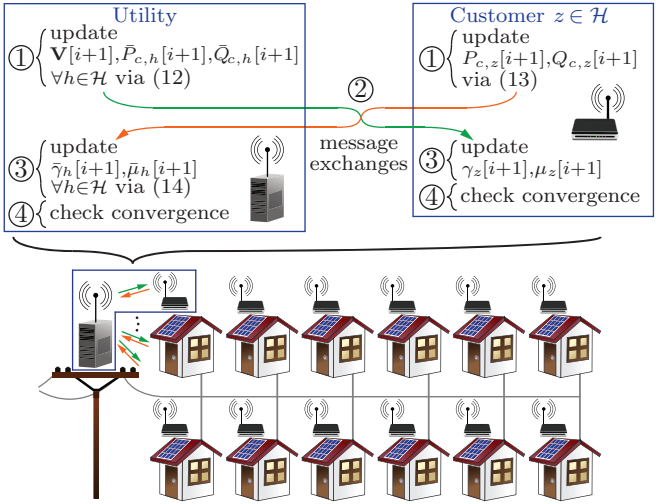}
\end{center}
\vspace{-0.25cm}
\caption{DOID: scenario with utility-customer
  message passing according to Algorithm 1.}
\label{Fig:DOID2}
\end{figure}

The resultant decentralized algorithm entails a two-way message exchange
between the utility and customers of the current iterates
$\bar{\bp}_{c}[i],\bar{\bq}_{c}[i]$ and
$\bp_{c}[i],\bq_{c}[i]$. Specifically, at each iteration $i > 0$, the
utility-owned device solves the OID rendition~\eqref{step1_utility} to
update the desired PV-inverter setpoints based on the performance
objectives described by $\bar{C}(\bV, \bar{\bp}_{c},\bar{\bq}_{c})$
(which is regularized with the term $F(\bar{\bp}_{c},\bar{\bq}_{c},
\{\cP_h[i]\})$ enforcing consensus with the setpoints computed at the
customer side), as well as the electrical network
constraints~\eqref{Pmc-balanceP}-\eqref{Pmc-Vlimits};
once~\eqref{step1_utility} is solved, the utility relays to each
customer a copy of the iterate value
$(\bar{P}_{c,h}[i+1], \bar{Q}_{c,h}[i+1])$. In the meantime, the
PV-inverter setpoints are simultaneously updated via~\eqref{step1_customer} and subsequently sent to the utility. Once
the updated local iterates are exchanged, utility and customers update
the local dual variables~\eqref{step2}.

The resultant decentralized algorithm is tabulated as Algorithm~1,
illustrated in Fig.~\ref{Fig:DOID2}, and
its convergence to the solution of the centralized OID
problem~\eqref{Poid} is formally stated next.

\begin{proposition}
\label{prop:convergence}
The iterates $\bar{\cP}[i], \{\cP_h[i]\}$ and $\cD[i]$ produced by
\textbf{[S1$^\prime$]}--\textbf{[S2$^\prime$]} are convergent, for any $\kappa >
0$. Further, $\lim_{i \rightarrow + \infty} \bV[i] =
\bV^{\mathrm{opt}}$, $\lim_{i \rightarrow + \infty} \bp_{c}[i] =
\lim_{i \rightarrow + \infty} \bar{\bp}_{c}[i] =
\bp_{c}^{\mathrm{opt}}$ and $\lim_{i \rightarrow + \infty} \bq_{c}[i]
= \lim_{i \rightarrow + \infty} \bar{\bq}_{c}[i] =
\bq_{c}^{\mathrm{opt}}$, with $\bV^{\mathrm{opt}},
\bp_{c}^{\mathrm{opt}}, \bq_{c}^{\mathrm{opt}}$ denoting the optimal solutions
of the OID problems~\eqref{Poid} and~\eqref{Pmconsensus1}.  \hfill
$\Box$
\end{proposition}

Notice that problem~\eqref{step1_utility} can be conveniently
reformulated in a standard SDP form (which involves the minimization
of a linear function, subject to linear (in)equalities and linear
matrix inequalities) by introducing pertinent auxiliary optimization
variables and by using the Schur
complement~\cite{Vandenberghe96,LavaeiLow,OID}. 
Finally, for a given consensus error $0 < \epsilon \ll 1$, the
algorithm terminates when $\|\bar{\bp}_{c}[i] - \bp_{c}[i]\|_2^2 +
\|\bar{\bq}_{c}[i] - \bq_{c}[i]\|_2^2 \leq \epsilon$. However, it is
worth emphasizing that, at each iteration $i$, the utility company
solves a consensus-enforcing regularized OID problem, which yields
intermediate voltages and power flows that clearly adhere to
electrical network constraints.

  Once the decentralized algorithm has converged, the
  real and reactive setpoints are implemented in the PV
  inverters. Notice however that Algorithm~1 affords an \emph{online}
  implementation; that is, the intermediate PV-inverter setpoints
  $\bar{\bp}_{c}[i], \bar{\bq}_{c}[i]$ are dispatched (and set at the
  customer side) as and when they become available, rather than
  waiting for the algorithm to converge.

\section{DOID: Network cluster partitions} \label{sec:Distributed1}

\begin{figure}[t]
\begin{center}
\includegraphics{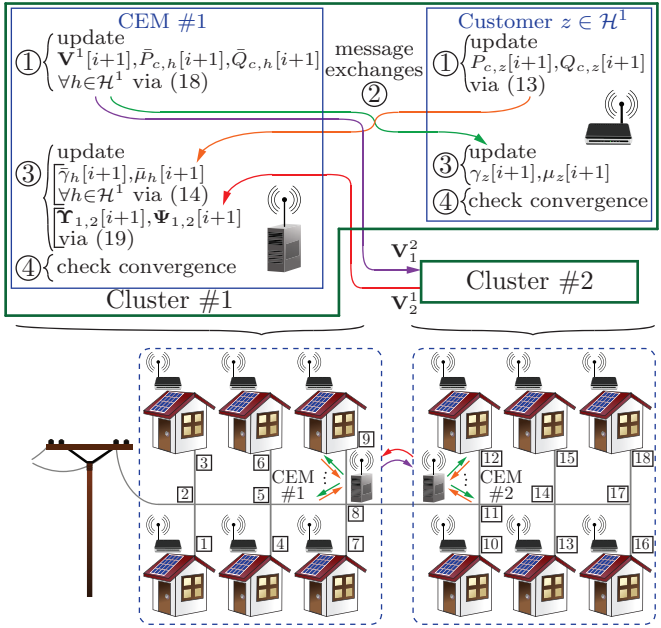}
\end{center}
\vspace{-0.25cm}
\caption{Network division into clusters and illustration of Algorithm 2. In this setup, $\cC^1 = \{1,\dots,9\}$, $\tilde\cC^1 = \{1,\dots,9,11\}$, $\cC^2 = \{10,\dots,18\}$, and $\tilde\cC^2 = \{8,10,\dots,18\}$.}
\label{Fig:DOID3}
\end{figure}

Consider the case where the distribution network is partitioned into
clusters, with $\cC^{a} \subset \cN$ denoting the set of nodes within
cluster $a$. Also, define $\tilde{\cC}^{a} := \cC^{a} \cup \{n | (m,n) \in \cE, m \in
\cC^{a}, n \in \cC^{j}, a \neq j \}$; that is, $\tilde{\cC}^{a}$
also includes the nodes belonging to different clusters that are
connected to the $a$-th one by a distribution
line~\cite{Zhu_DSE,Dallanese-TSG13} (see Fig.~\ref{Fig:DOID3} for an illustration). Hereafter, superscript
$(\cdot)^{a}$ will be used to specify quantities pertaining to cluster
$a$; e.g., $\cH^{a}$ is the set of houses located within
cluster $\cC^{a}$, and vectors $\bar{\bp}_{c}^{a},
\bar{\bq}_{c}^{a}$ collect copies of the setpoints of PV inverters $h
\in \cH^a$ available with the $a$-th CEM [cf.~\eqref{Pmconsensus1}]. With regard to notation, an
exception is $\bV^{a}$, which denotes the sub-matrix of
$\bV$ corresponding to nodes in the extended cluster
$\tilde{\cC}^{a}$.

Based on this network partitioning, consider decoupling the
network-related cost $\bar{C}(\bV, \bar{\bp}_{c},\bar{\bq}_{c})$
in~\eqref{Pmc-cost} as
\begin{align}
\bar{C}(\bV, \bar{\bp}_{c},\bar{\bq}_{c}) = \sum_{a = 1}^{N_a} \underbrace{ \left[C^a(\bV^a, \bar{\bp}_{c}^a) + \lambda^{a} \sum_{h \in \cH^{a}} \|[\bar{P}_{c,h}, \bar{Q}_{c,h}]\|_2 \right]}_{:= \bar{C}^a(\bV^a, \bar{\bp}_{c}^a,\bar{\bq}_{c}^a)} \nonumber 
\end{align} 
where $N_a$ is the number of clusters, $C^a(\bV^a, \bar{\bp}_{c}^a)$ captures optimization objectives
of the $a$-th cluster (e.g., power losses within the
cluster~\cite{Baldick99,Dallanese-TSG13}), and the sparsity-promoting
regularization function is used to determine which PV inverters in
$\cH^a$ provide ancillary services. Further, per-cluster $a = 1, \ldots, N_a$,
define the region of feasible power flows as
[cf.~\eqref{Pmc-balanceP}-\eqref{Pmc-Vlimits}]:
\begin{align}
\cR^{a}  :=  
\left\{ \bV^{a}, \bar{\bp}_{c}^{a}, \bar{\bq}_{c}^{a}:  \hspace{-.2cm} 
\begin{array}{l}
\trace(\bA_h^{a} \bV^{a}) = \bar{P}_{h} - \bar{P}_{c,h} , \forall h \in \cH^{a}   \\
\trace(\bB_h^{a} \bV^{a})  = - Q_{\ell,h} + \bar{Q}_{c,h} , \forall h \in \cH^{a}    \\
\trace(\bA_n^{a}  \bV^{a} )  = 0,  \forall \, n \in \cU^{a}  \\
\trace(\bB_n^{a}  \bV^{a} ) = 0, \forall \, n \in \cU^{a} \\
V_{\mathrm{min}}^2  \leq \trace(\bM_n^{a} \bV^{a})  \leq V_{\mathrm{max}}^2, \forall \, n \in \cC^{a} 
\end{array}
\hspace{-.2cm} \right\} \nonumber 
\end{align}
where $\bA_h^{a}, \bB_h^{a}$, and $\bM_h^{a}$ are the sub-matrices of
$\bA_h, \bB_h$, and $\bM_h$, respectively, formed by extracting rows
and columns corresponding to nodes in $\tilde{\cC}^{a}$. With these
definitions, problem~\eqref{Pmconsensus1} can be equivalently
formulated as:
 \begin{subequations} 
\label{PmAreaConsensus1}
\begin{align}
& \hspace{-2.2cm} \min_{\substack{\bV, \bp_{c}, \bq_{c} \\ \bar{\bp}_{c},\bar{\bq}_{c}}} \,\,  \sum_{a} \left[\bar{C}^a(\bV^a, \bar{\bp}_{c}^a,\bar{\bq}_{c}^a) + \sum_{h \in \cH^a} R_{h}(P_{c,h}) \right] \label{Pmca-cost} \\
\,\, \mathrm{s.\,to} \,\,  \bV  \succeq \mathbf{0}, &  \mathrm{~and}  \nonumber  \\ 
& \hspace{-.8cm} \{\bV^{a}, \bar{\bp}_{c}^{a},  \bar{\bq}_{c}^{a}\}  \in \cR^{a}  \hspace{1.55cm} \forall \,\, a \label{Pmca-region}  \\
& \hspace{-.8cm} (P_{c,h}, Q_{c,h}) \in \cF^\mathrm{OID}_h  \hspace{1.3cm} \forall \, h \in \cH^{a}, \forall \,\, a  \label{Pmca-oid}  \\
\bar{P}_{c,h} & = P_{c,h}, \,\, \bar{Q}_{c,h}  =  Q_{c,h} \hspace{.4cm} \forall \,  h \in \cH^{a}, \forall \,\, a .   \label{Pmca-consensusP} 
\end{align} 
\end{subequations}
Notice that, similar to~\eqref{Pmc-consensusP},
constraints~\eqref{Pmca-consensusP} ensure that the CEM and
customer-owned PV systems consent on the optimal PV-inverter
setpoints. Formulation~\eqref{PmAreaConsensus1} effectively decouples
cost, power flow constraints, and PV-related consensus
constraints~\eqref{Pmca-consensusP} on a per-cluster basis. The main challenge towards solving~\eqref{PmAreaConsensus1} in a
decentralized fashion lies in the positive semidefinite (PSD) constraint $\bV \succeq \mathbf{0}$, which
clearly couples the matrices $\{\bV^{a}\}$. To address this challenge, results on completing partial Hermitian matrices from~\cite{Grone84} will be leveraged to identify partitions of the distribution
network in clusters for which the PSD constraint on $\bV$ would
decouple to $\bV^{a} \succeq \mathbf{0}$, $\forall a$. This decoupling would clearly facilitate
the decomposability of~\eqref{PmAreaConsensus1} in per-cluster
sub-problems~\cite{Zhu_DSE,Dallanese-TSG13}.

Towards this end, first define the set of \emph{neighboring clusters}  for
the $a$-th one as $\tilde{\cB}^{a} := \{j | \tilde{\cC}^{a} \cap
\tilde{\cC}^{j}\neq 0\}$. Further, let $\cG_{\cC}$ be a graph
capturing the control architecture of the distribution network, where
nodes represent the clusters and edges connect neighboring clusters
(i.e., based on sets $\{\tilde{\cB}^{a}\}$); for example, the graph  $\cG_{\cC}$ associated with the network in 
Fig.~\ref{Fig:DOID3} has two nodes, connected through an edge (since the two areas are connected). In general, it is clear that if clusters $a$
and $j$ are neighbors, then CEM $a$ and CEM $j$ must agree on the
voltages at the two end points of the distribution line connecting the
two clusters. For example, with reference to Fig.~\ref{Fig:DOID3}, notice that line $(8,11)$ connects clusters $1$ and $2$. Therefore, CEM 1 and CEM 2 must agree on voltages $V_8$ and $V_{11}$. Lastly,
let $\bV_{j}^{a}$ denote the sub-matrix of $\bV^{a}$ corresponding to
the two voltages on the line connecting clusters $a$ and
$j$. Recalling the previous example, agreeing on $V_8$ and $V_{11}$ is
tantamount to setting $\bV_{2}^{1} = \bV_{1}^{2}$, where $\bV_{2}^{1}$
is a $2 \times 2$ matrix representing the outer-product
$[V_8,V_{11}]^\sfT[V_8,V_{11}]^*$. Using these definitions, the following proposition can be proved by
suitably adapting the results of~\cite{Zhu_DSE,Dallanese-TSG13} to the
problem at hand.
\begin{proposition}
\label{prop:equivalentSDP}
Suppose: (i) the cluster graph $\cG_{\cC}$ is a tree, and (ii)
clusters are not nested (i.e., $|\tilde{\cC}^{a} \backslash
(\tilde{\cC}^{a} \bigcap \tilde{\cC}^{j} ) | > 0$ $\forall a \neq
j$). Then,~\eqref{PmAreaConsensus1} is equivalent to the following
 problem:
 \begin{subequations} 
\label{PmAreaConsensus2}
\begin{align}
& \hspace{-1.2cm} \min_{\substack{\{\bV^a, \bp_{c}^a, \bq_{c}^a\} \\ \bar{\bp}_{c},\bar{\bq}_{c}}} \,\,  \sum_{a} \left[\bar{C}^a(\bV^a, \bar{\bp}_{c}^a,\bar{\bq}_{c}^a) + \sum_{h \in \cH^a} R_{h}(P_{c,h}) \right] \label{Pmca-cost2} \\
\,\, \mathrm{s.\,to} & \,\, \eqref{Pmca-region}-\eqref{Pmca-consensusP}  \mathrm{~and}  \nonumber  \\ 
& \bV^{a} \succeq \mathbf{0}  \hspace{.85cm} \forall \,\, a \label{localPSD} \\
&  \bV_{j}^{a} = \bV^{j}_a, \hspace{.5cm} \forall \,  j \in \tilde{\cB}^{a}, \,\, \forall \, a.  \label{consensusVoltages} 
\end{align} 
\end{subequations}
Under (i)--(ii), there exists a rank-$1$ matrix $\bV^\textrm{opt}$
solving~\eqref{PmAreaConsensus1} optimally if and only if
$\rank\{\bV^a\} = 1,\,\forall a = 1, \ldots, N_a$.  \hfill
$\Box$
\end{proposition}

Notice that the $|\cN| \times |\cN|$ matrix $\bV$ is replaced by per-cluster reduced-dimensional $|\tilde{\cC}^{a}| \times |\tilde{\cC}^{a}|$ matrices $\{\bV^a\}$ in~\eqref{PmAreaConsensus2}.
Proposition~\ref{prop:equivalentSDP} is grounded on the results
of~\cite{Grone84}, which asserts that a PSD matrix $\bV$ can be
obtained starting from sub-matrices $\{\bV^a\}$ if and only if the
graph induced by $\{\bV^a\}$ is chordal. Since a PSD matrix can be reconstructed from 
$\{\bV^a\}$, it suffices to impose contraints $\bV^{a} \succeq \mathbf{0}$,
$\forall a = 1, \ldots, N_a$. Assumptions
(i)--(ii) provide sufficient conditions for the graph induced by
$\{\bV^a\}$ to be chordal, and they are typically satisfied in practice (e.g., when each cluster is set to be a lateral or a sub-lateral). The second part of the proposition asserts
that, for the completable PSD matrix $\bV$ to have rank $1$, all
matrices $\bV^{a}$ must have rank $1$; thus, if $\rank\{\bV^a\} = 1$
for all clusters, then $\{\bV^a\}$ represents a globally
optimal power flow solution for given inverter setpoints. 

Similar to~\eqref{Pmconsensus2}, auxiliary variables are introduced to
enable decomposability of~\eqref{PmAreaConsensus2} in per-cluster
subproblems. With variables $x_h, y_h$ associated with inverter $h$,
and $\bW^{a,j}, \bQ^{a,j}$ with neighboring clusters
$a$ and $j$,~\eqref{PmAreaConsensus2} is reformulated as:
 \begin{subequations} 
\label{PmAreaConsensus3}
\begin{align}
& \hspace{-1.9cm} \min_{\substack{\{\bV^a, \bp_{c}^a, \bq_{c}^a\} \\ \bar{\bp}_{c}^a,\bar{\bq}_{c}^a\\
\{\bW^{a,j}, \bQ^{a,j} , x_h, y_h\}}} \,\,  \sum_{a} \left[\bar{C}^a(\bV^a, \bar{\bp}_{c}^a,\bar{\bq}_{c}^a) + \sum_{h \in \cH^a} R_{h}(P_{c,h}) \right] \nonumber \\
\,\, \mathrm{s.\,to} \,\, \eqref{Pmca-region}&-\eqref{Pmca-oid} ,  \bV^{a} \succeq \mathbf{0}  \,\, \forall a, \mathrm{~and}  \nonumber  \\ 
 \Re\{\bV^{a}_{j}\} & = \bW^{a,j},  \,\,   \bW^{a,j} = \bW^{j,a} \hspace{.3cm} \forall \,  j \in \tilde{\cB}^{a}, \,\, \forall \,  a  \label{PmAreaConsensus3-real} \\
  \Im\{\bV^{a}_{j}\} & = \bQ^{a,j}, \hspace{.35cm}  \bQ^{a,j} = \bQ^{j,a} \hspace{.4cm} \forall \,  j \in \tilde{\cB}^{a}, \,\, \forall \, a \label{PmAreaConsensus3-imag} \\
\bar{P}_{c,h} & = x_h, \hspace{.7cm}  x_h = P_{c,h} \hspace{.75cm} \forall \, h \in \cH^{a}, \,\, \forall \, a   \label{Pmc-consP2} \\
\bar{Q}_{c,h} & =  y_h, \hspace{.75cm}  y_h =  Q_{c,h}  \hspace{.75cm} \forall \, h \in \cH^{a}, \,\, \forall \, a.  \label{Pmc-consQ2} 
\end{align} 
\end{subequations}
This problem can be solved across clusters by resorting
to ADMM. To this end, a partial quadratically-augmented
Lagrangian, obtained by dualizing constraints $\Re\{\bV^{a}_{j}\} =
\bW^{a,j}$, $\Im\{\bV^{a}_{j}\} = \bQ^{a,j}$,
$\bar{P}_{c,h} = x_h$, and $\bar{Q}_{c,h} = y_h$ is defined first;
then, the standard ADMM steps involve a cyclic minimization of the
resultant Lagrangian with respect to $\{\bV^a, \bp_{c}^a, \bq_{c}^a,
\bar{\bp}_{c}^a,\bar{\bq}_{c}^a\}$ (by keeping the remaining variables fixed); the auxiliary variables $\{\bW^{a,j},
\bQ^{a,j} , x_h, y_h\}$; and, finally, a dual ascent
step~\cite[Sec.~3.4]{BeT89}. It turns out that Lemma~\ref{lemma:dualupdates} still holds in the
present case. Thus, using this lemma, along with the result
in~\cite[Lemma~3]{Dallanese-TSG13}, it can be shown that the ADMM
steps can be simplified as described next (the derivation is omitted due to space limitations):

\noindent \textbf{[S1$^{\prime \prime}$]}  Each PV system updates the local copy $\cP_h[i+1]$ via~\eqref{step1_customer}; while, each CEM updates the voltage profile of its cluster, and the local copies of the setpoints of inverters $\cH^{a}$ by solving the following convex problem:
\begin{subequations} 
\label{Cluster_step1} 
\begin{align}
& \hspace{-.1cm}  \bar{\cP}^{a}[i+1] := \arg  \min_{\substack{\bV^{a}, \bar{\bp}_{c}^{a}, \bar{\bq}_{c}^{a} \\ \{\alpha_j \geq 0, \beta_j \geq 0\}}} \Big[
\bar{C}^a(\bV^a, \bar{\bp}_{c}^a,\bar{\bq}_{c}^a)   \nonumber \\ 
&  \hspace{1.3cm} + F^{a}(\bar{\bp}_{c}^{a},\bar{\bq}_{c}^{a}, \{\cP_h[i]\}) + F_V^{a}(\bV^{a}, \{\bV^{j}[i]\})  \Big]\label{step1_utility_2} \\
& \hspace{1.2cm} \mathrm{s.\,to~} \{\bV^{a}, \bar{\bp}_{c}^{a},  \bar{\bq}_{c}^{a}\}  \in \cR^{a}, \bV^{a} \succeq \mathbf{0}, \mathrm{~and:}  \\
& \hspace{1.8cm} 
\left[ \begin{array}{cc}
-\alpha_{j} & \ba_{j}^\sfT \\
\ba_{j} & -\bI
\end{array} \right] \preceq \mathbf{0}, \,\,\, \forall j \in \tilde{\cB}^{a} \label{penalty} \\
& \hspace{1.8cm} 
\left[ \begin{array}{cc}
-\beta_{j} & \bb_{j}^\sfT \\
\bb_{j}  & -\bI
\end{array} \right] \preceq \mathbf{0}, \,\,\, \forall j \in  \tilde{\cB}^{a} \label{penalty} 
\end{align}
where vectors $\ba_{j}$ and $\bb_{j}$ collect the real and imaginary parts, respectively, of the entries of the matrix $\bV^{a}_j - \frac{1}{2}\left(\bV^{a}_j[i] + \bV^{j}_a[i] \right)$;  the regularization function $F^{a}(\bar{\bp}_{c}^{a},\bar{\bq}_{c}^{a}, \{\cP_h[i]\})$ enforcing consensus on the inverter setpoints is defined as in~\eqref{step1_utility_3} (but with the summation limited to inverters $\cH^{a}$) and,  $F_V(\bV^{a}, \{\bV^{j}[i]\})$ is given by:
\begin{align}
F_V^{a}(\bV^{a}, \{\bV^{j}[i]\})& := \sum_{j \in \tilde{\cB}^{a}} \Big[ \frac{\kappa}{2}(\alpha_{j}  + \beta_{j} ) +
\trace( \bUpsilon_{a,i}^\sfT[i] \Re\{\bV^{a}_{j}\}) \nonumber \\
& \hspace{1.2cm} + \trace( \bPsi_{a,i}^\sfT[i] \Im\{\bV^{a}_{j}\})  \Big] \, .
\end{align}
\end{subequations}

\noindent \textbf{[S2$^{\prime \prime}$]} Update dual variables $\{\gamma_h, \mu_h\}$ via~\eqref{step2_customer_and_utility} at both, the customer and the CEMs; variables $\{ \bUpsilon_{{a,i}},  \bPsi_{{a,i}}\}$ are updated locally per cluster $a = 1,\ldots,N_a$ as:
\small
\begin{subequations} 
\label{Cluster_step2} 
\begin{align}
&\hspace{-.2cm}  \bUpsilon_{{a,j}}[i+1] = \bUpsilon_{{a,j}}[i] + \frac{\kappa}{2} \left( \Re\{\bV^{a}_{j}[i+1]\} -  \Re\{\bV^{j}_{a}[i+1]\}\right)  \\
&\hspace{-.2cm}  \bPsi_{{a,j}}[i+1] = \bPsi_{{a,j}}[i] + \frac{\kappa}{2} \left(  \Im\{\bV^{a}_{j}[i+1]\} - \Im\{\bV^{j}_{a}[i+1]\}  \right) . 
\end{align}
\end{subequations}
\normalsize

\begin{algorithm}[t!]
\label{alg:OID3}
\caption{DOID: multi-cluster distributed optimization} \small{
\begin{algorithmic}

\STATE Set $\gamma_h[0]  = \mu_h[0] = 0$ for all $h \in \cH^a$ and for all clusters.

\STATE Set $\bUpsilon_{{a,j}}[0]  =  \bPsi_{{a,j}}[0] = \mathbf{0}$  for all pair of neighboring clusters.

\FOR {$i = 1,2,\ldots$ (repeat until convergence)} 

\STATE 1. \textbf{[CEM-$a$]}: update $\bV^a[i+1]$ and $ \bar{\bp}_{c}^{a}, \bar{\bq}_{c}^{a}$ via~\eqref{Cluster_step1}. \\

\hspace{.25cm}  \textbf{[Customer-$h$]}: update $\bar{P}_{c,h}[i+1], \bar{Q}_{c,h}[i+1]$ via~\eqref{step1_customer}. \\

\STATE 2.  \textbf{[CEM-$a$]}: send $\bV^{a}_{j}[i+1]$ to CEM $j$;

\hspace{.25cm} \textbf{[CEM-$a$]}: receive $\bV^{j}_{a}[i+1]$ from CEM $j$;

\hspace{1.7cm} repeat $\forall j \in \tilde{\cB}^{a}$

\STATE 3.  \textbf{[CEM-$a$]}: send  $\bar{P}_{c,h}[i+1], \bar{Q}_{c,h}[i+1]$ to EMU-$h$; 

\hspace{1.5cm} repeat for all $h \in \cH^a$.

\hspace{.2cm}  \textbf{[Customer-$h$]}: receive $\bar{P}_{c,h}[i+1], \bar{Q}_{c,h}[i+1]$ from CEM $a$; 

\hspace{1.65cm} send $P_{c,h}[i+1], Q_{c,h}[i+1]$ to CEM $a$;

\hspace{1.65cm} repeat for all $h \in \cH^a$.

\hspace{.2cm}  \textbf{[CEM-$a$]}: receive $P_{c,h}[i+1], Q_{c,h}[i+1]$ from $h$;  

\hspace{1.5cm} repeat for all $h \in \cH^a$.

\STATE 4. \textbf{[CEM-$a$]}: update $\{\gamma_h[i+1] , \mu_h[i+1] \}_{h \in \cH}$ via~\eqref{step2}.

\hspace{.2cm} \textbf{[CEM-$a$]}: update $\{\bUpsilon_{{a,j}}[i+1], \bPsi_{{a,j}}[i+1]\}$ via~\eqref{Cluster_step2};

\hspace{.2cm}  \textbf{[Customer-$h$]}: update dual variables $\gamma_h[i+1] , \mu_h[i+1]$ via~\eqref{step2};

\ENDFOR

Implement setpoints in the PV inverters. 

\end{algorithmic}}
\end{algorithm}

The resultant decentralized algorithm is tabulated as Algorithm~2,
illustrated in Fig. \ref{Fig:DOID3}, and it involves an exchange of:
(i) the local submatrices $\{\bV^{a}_{j}[i+1]\}$ among neighboring
CEMs to agree upon the voltages on lines connecting clusters; and,
(ii) the local copies of the PV inverter setpoints between the CEM and
customer-owned PV systems. Using arguments similar to
Proposition~\ref{prop:convergence}, convergence of the algorithm can
be readily established.

\begin{proposition}
\label{prop:convergence}
For any $\kappa > 0$, the iterates $\{\bar{\cP}^{a}[i]\},
\{\cP_h[i]\}, \cD[i]$ produced by \textbf{[S1$^{\prime \prime}$]}--\textbf{[S2$^{\prime \prime}$]} are convergent, and they converge to a solution of the OID
problems~\eqref{Poid} and~\eqref{PmAreaConsensus1}.  \hfill $\Box$
\end{proposition}

Once the decentralized algorithm has converged, the real and reactive setpoints are implemented by the PV inverter controllers.

Finally, notice that the worst case complexity of an SDP is on the order $\mathcal{O}(\max\{N_c,N_v\}^4 \sqrt{N_v} \log(1/\epsilon))$ for general purpose solvers, with $N_c$ denoting the total number of constraints, $N_v$ the total number of variables, and $\epsilon > 0$ a given solution accuracy~\cite{Vandenberghe96}. It follows that the worst case complexity of~\eqref{Cluster_step1} is markedly lower than the one of the centralized problem~\eqref{Poid}. Further, the sparsity of $\{\bA_{n},\bB_{n}, \bM_{n}\}$ and the so-called chordal structure of the underlying electrical graph matrix can be exploited to obtain substantial computational savings; see e.g.,~\cite{Jabr12}.

\section{Case Studies} \label{sec:Simulations}

Consider the distribution network in Fig.~\ref{F_LV_network}, which is adopted
from~\cite{Tonkoski11,OID}. The simulation
parameters are set as in~\cite{OID} to check the consistency
between the results of centralized and decentralized
schemes. Specifically, the pole-pole distance is set to $50$ m;
lengths of the drop lines are set to $20$ m; and voltage limits
$V^\textrm{min}, V^\textrm{max}$ are set to 0.917 pu and 1.042 pu,
respectively (see e.g.,~\cite{Tonkoski11}).  The
optimization package \texttt{CVX}\footnote{[Online] Available:
  \texttt{http://cvxr.com/cvx/}} is employed to solve
relevant optimization problems in \texttt{MATLAB}. In all the conducted
numerical tests, the rank of matrices $\bV$ and $\{\bV^{a}\}$ was
always $1$, meaning that globally optimal power flow solutions were
obtained for given inverter setpoints.

\begin{figure*}[t]
\begin{center}
\subfigure[]{\includegraphics[width=8.5cm]{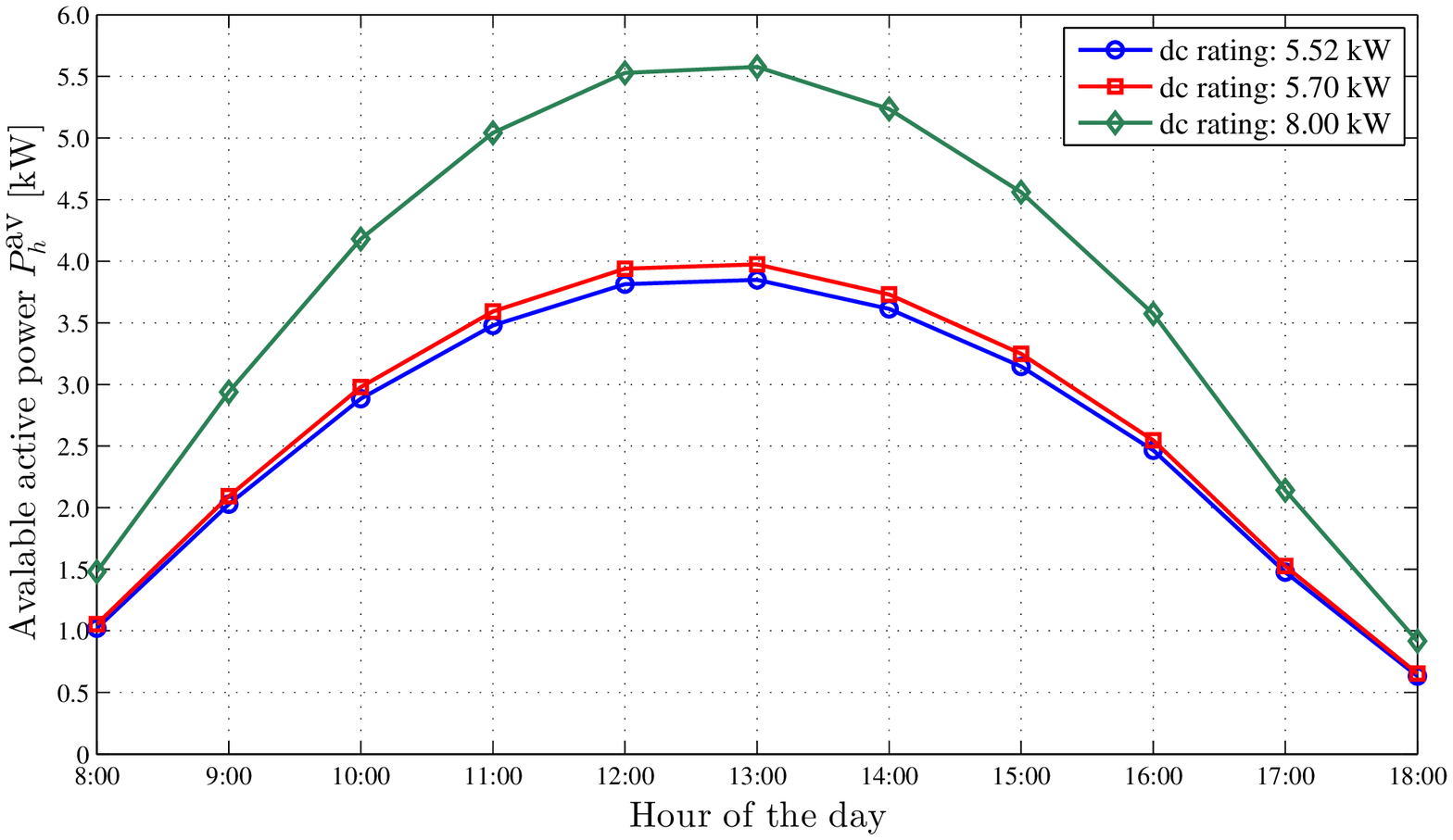}}
\subfigure[]{\includegraphics[width=8.5cm]{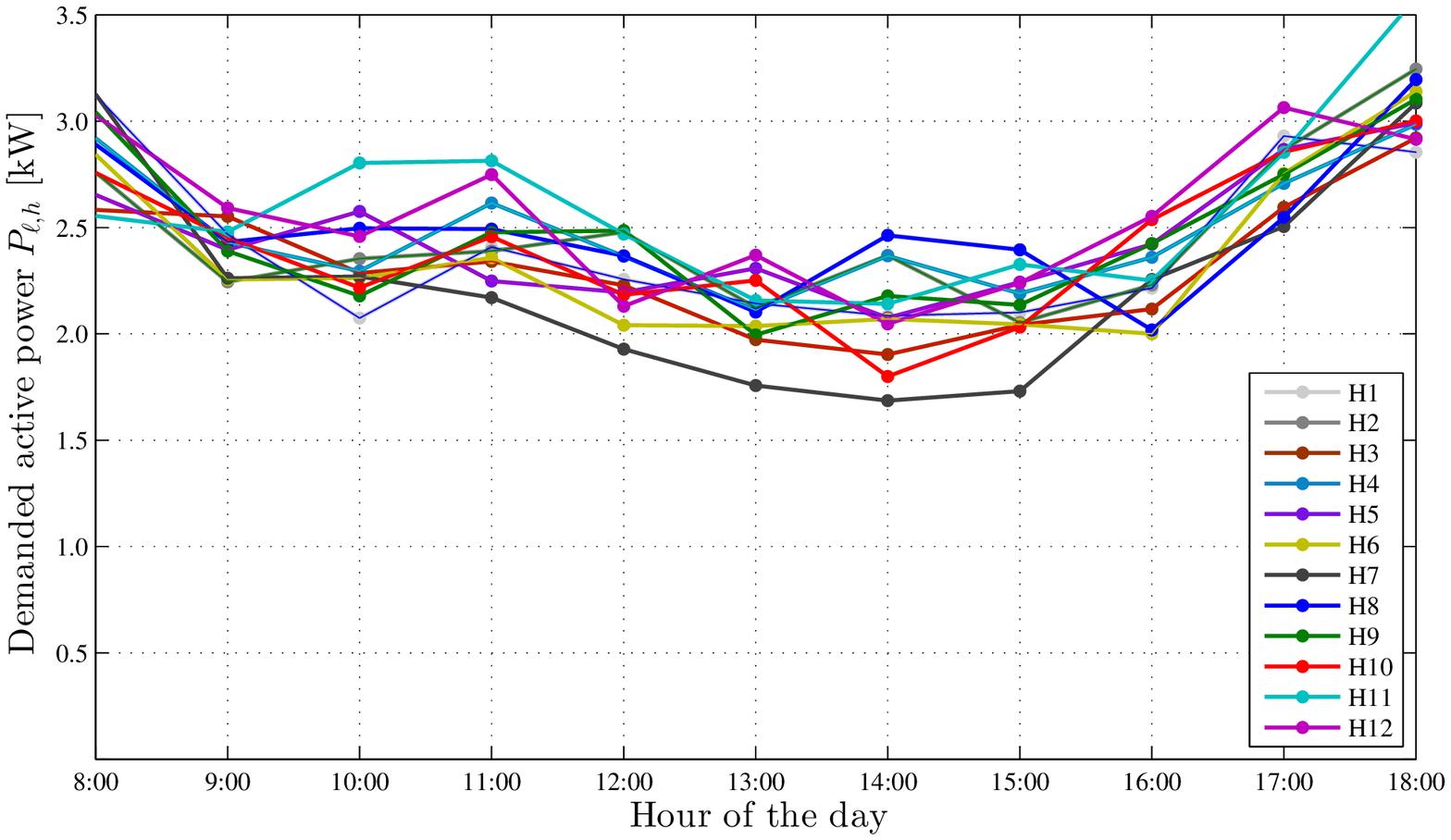}} 
\end{center}
\vspace{-.3cm}
\caption{Problem inputs: (a) available active powers $\{P_h^\textrm{av}\}$ from inverters with dc ratings of $5.52$ kW, $5.70$ kW, and $8.00$ kW; (b) demanded active loads at the households (reactive demand is computed by presuming a power factor of 0.9). }
\vspace{-.3cm}
\label{Fig:OIDPowers}
\end{figure*}

The available
active powers $\{P_h^\textrm{av}\}_{h \in \cH}$ are computed using the System
Advisor Model (SAM)\footnote{[Online] Available at
  \texttt{https://sam.nrel.gov/}.} of the National Renewable Energy
Laboratory (NREL); specifically, the typical meteorological year (TMY)
data for Minneapolis, MN, during the month of July are used.
All 12 houses feature fixed roof-top PV systems, with a dc-ac derating
coefficient of $0.77$. The dc ratings of the houses are as follows:
$5.52$ kW for houses $\mathrm{H}_1, \mathrm{H}_9, \mathrm{H}_{10}$;
$5.70$ kW for $\mathrm{H}_2, \mathrm{H}_6, \mathrm{H}_{8},
\mathrm{H}_{11}$; and, $8.00$ kW for the remaining five houses. 
The active powers $\{P_h^\textrm{av}\}$ generated by 
the inverters with dc ratings of $5.52$ kW, $5.70$ kW, and $8.00$ kW
are plotted in Fig.~\ref{Fig:OIDPowers}(a). As
suggested in~\cite{Turitsyn11}, it is assumed that the PV inverters
are oversized by $10\%$ of the resultant ac rating. The minimum power
factor for the inverters is set to 0.85~\cite{Braun10}.

\begin{figure*}[t]
\begin{center}
\subfigure[]{\includegraphics[width=8.5cm]{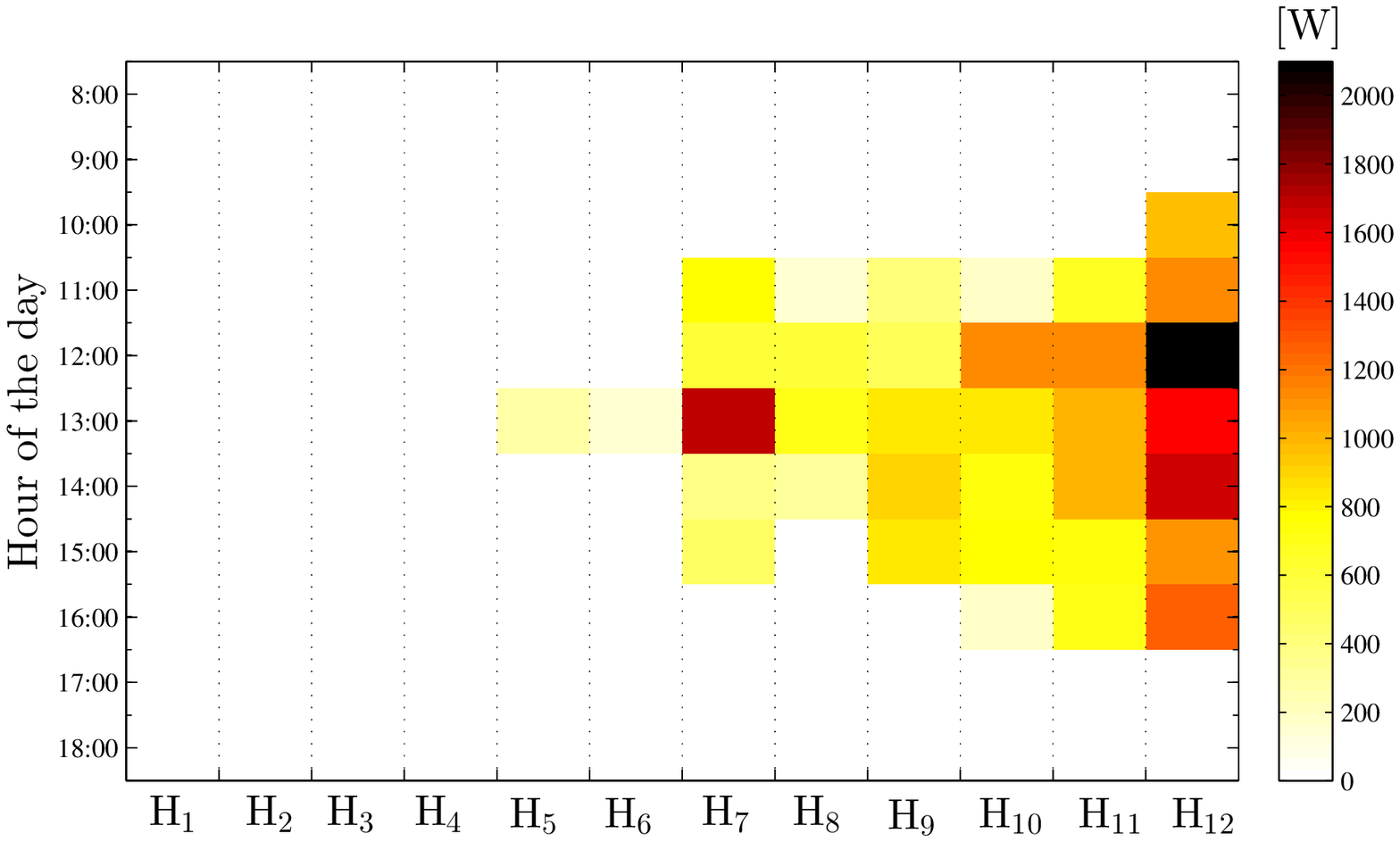}}
\subfigure[]{\includegraphics[width=8.5cm]{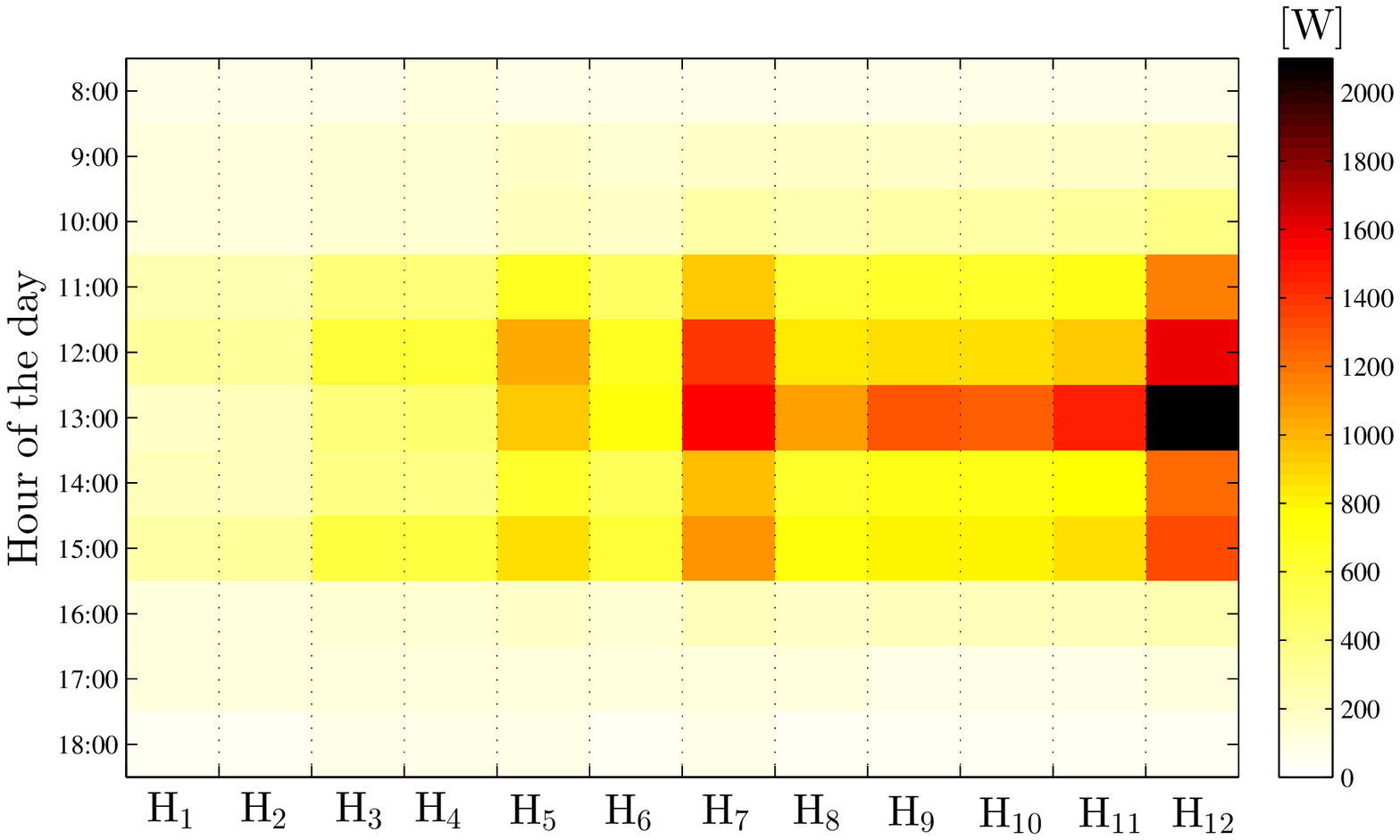}} 
\end{center}
\vspace{-.3cm}
\caption{Solution of the centralized OID problem~\eqref{Poid}: Curtailed active power per each household, for (a) $\lambda = 0.8$ and (b) $\lambda = 0$; see also~\cite{OID}.}
\label{Fig:OIDCentralizedResults}
\end{figure*}

\begin{figure*}[t]
\begin{center}
\subfigure[]{\includegraphics[width=8.5cm]{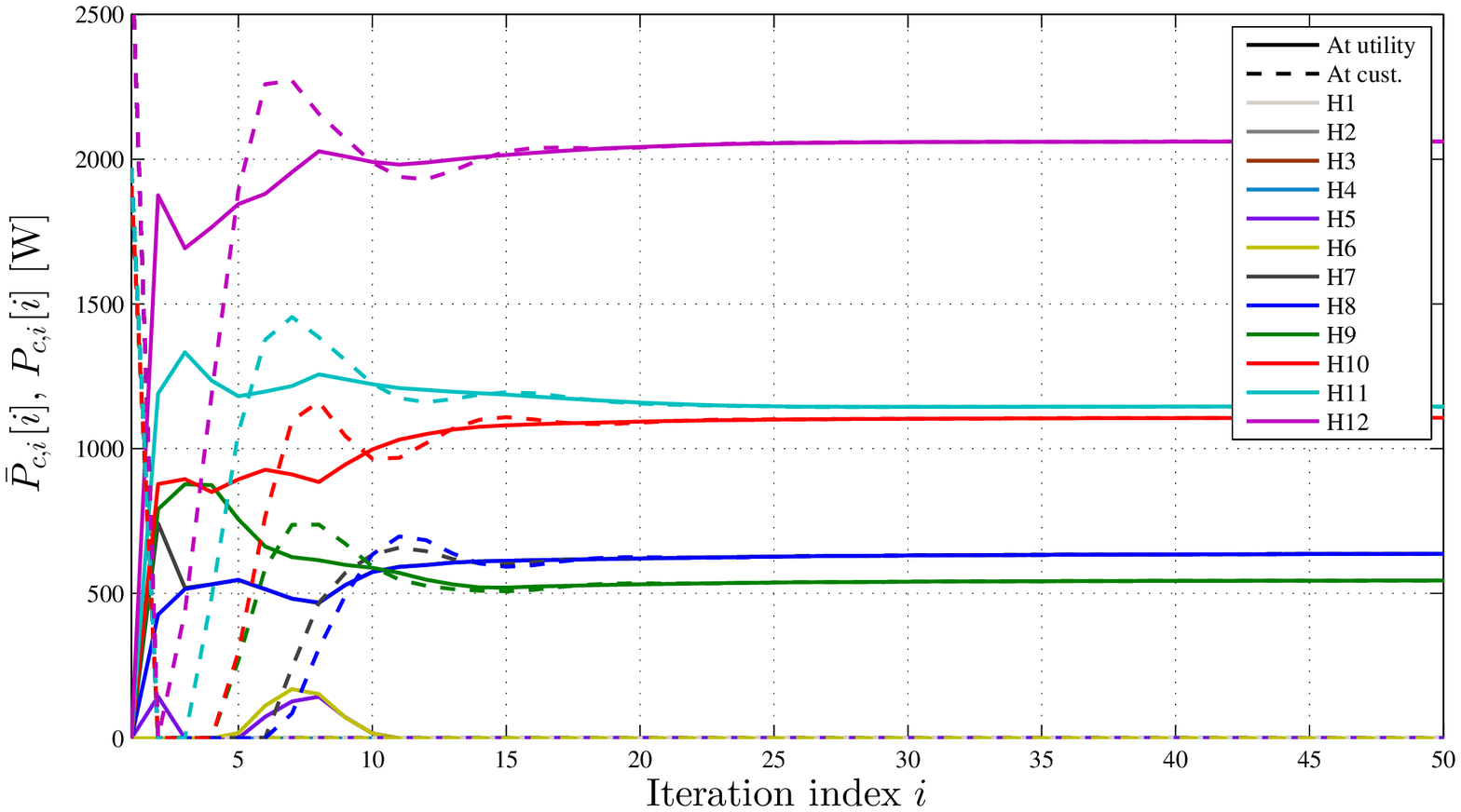}} 
\subfigure[]{\includegraphics[width=8.5cm]{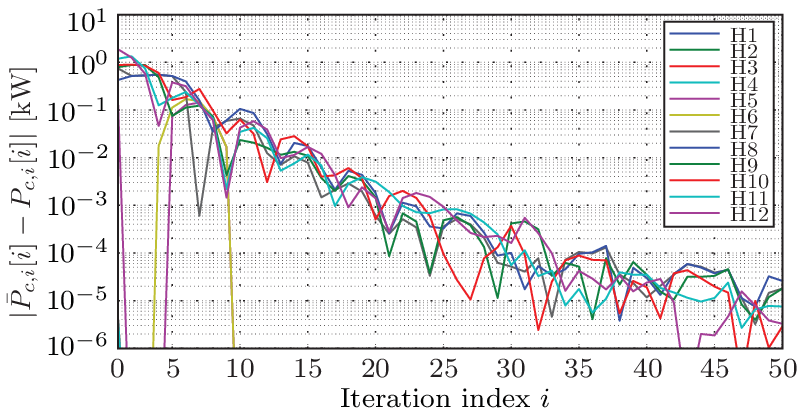}}
\end{center}
\vspace{-.3cm}
\caption{Convergence of Algorithm~1: (a) values of $\{P_{c,h}[i]\}_{h \in \cH}$ (dashed lines) and $\{\bar{P}_{c,h}[i]\}_{h \in \cH}$ as a function of the ADMM iteration index $i$. (b) Consensus error $|P_{c,h}[i] - \bar{P}_{c,h}[i]|$, for all houses $h \in \cH$ as a function of $i$.   }
\label{Fig:OID2results}
\end{figure*}

The residential load profile is obtained from the Open Energy Info
database and the base load experienced in downtown Saint Paul, MN,
during the month of July is used for this test case. To generate 12 different load profiles, 
the base active power profile is perturbed using a Gaussian random variable 
with zero mean and standard deviation $200$ W; the resultant active loads $\{P_{\ell,h}\}$
are plotted in Fig.~\ref{Fig:OIDPowers}(b). To compute the reactive loads $\{Q_{\ell,h}\}$, 
a power factor of 0.9 is presumed~\cite{Tonkoski11}.

Assume that the objective of the utility company is to minimize the power losses in the network; that is, upon defining the symmetric matrix $\bL_{mn} := \Re\{y_{mn}\} (\be_m - \be_n) (\be_m - \be_n)^\sfT$ per  distribution line $(m,n) \in \cE$, function $C_{\textrm{utility}}(\bV, \bar{\bp}_{c})$ is set to $ C_{\textrm{utility}}(\bV, \bar{\bp}_{c}) =  \trace(\bL \bV)$, with $\bL := \sum_{(m,n) \in \cE} \bL_{mn}$ (see~\cite{OID} for more details). At the customer side, function $R_{h}(P_{c,h})$ is set to $R_{h}(P_{c,h}) = 0.1 P_{c,h}$. The impact of varying the parameter $\lambda$ is investigated in detail in~\cite{OID}, and further illustrated in Fig.~\ref{Fig:OIDCentralizedResults}, where the solution of the centralized OID problem~\eqref{Poid} is reported for different values of the parameter $\lambda$ [cf.~\eqref{Glasso_powers}]. Specifically, Fig.~\ref{Fig:OIDCentralizedResults}(a) illustrates the active power curtailed from each inverter during the course of the day when $\lambda = 0.8$,  whereas the result in Fig.~\ref{Fig:OIDCentralizedResults}(b) were obtained by setting $\lambda = 0$. It is clearly seen that in the second case all inverters are controlled; in fact, they all curtail active power from 8:00 to 18:00. When $\lambda = 0.8$, the OID seeks a trade off between achievable objective and number of controlled inverters. It is clearly seen that the number of participating inverters grows with increasing solar irradiation, with a maximum of $7$ inverters operating away from the business-as-usual point at 13:00. 

The convergence of Algorithm~1 is showcased for $\lambda = 0.8$,
$C_{\textrm{utility}}(\bV, \bar{\bp}_{c}) = \trace(\bL \bV)$, and
$R_{h}(P_{c,h}) = 0.1 P_{c,h}$, and by utilizing the solar irradiation
conditions at 12:00. Figure~\ref{Fig:OID2results}(a) depicts the
trajectories of the iterates $\{P_{c,h}[i]\}_{h \in \cH}$ (dashed
lines) and $\{\bar{P}_{c,h}[i]\}_{h \in \cH}$ (solid lines) for all
the houses $\mathrm{H}_1-\mathrm{H}_{12}$. The
  results match the ones in~Fig.~\ref{Fig:OIDCentralizedResults}(a);
  in fact, at convergence (i.e., for iterations $i \geq 20$), only
  inverters at houses $\mathrm{H}_7-\mathrm{H}_{12}$ are controlled,
  and the active power curtailment set-points are in agreement.  This
  result is expected, since problems~\eqref{Poid}
  and~\eqref{Pmconsensus1} are \emph{equivalent}; the only difference
  is that~\eqref{Poid} affords only in centralized solution,
  whereas~\eqref{Pmconsensus1} is in a form that is suitable for the
  application of the ADMM to derive distributed solution schemes; see
  also~\cite{BoydADMoM,ErsegheADMM,ZhuADMoM}. Finally, the
trajectories of the set-point consensus error $|P_{c,h}[i] -
\bar{P}_{c,h}[i]|$, as a function of the ADMM iteration index $i$ are
depicted in Fig.~\ref{Fig:OID2results}(b). It can be clearly seen that
the algorithm converges fast to a set-point that is convenient for
both utility and customers. Similar trajectories were obtained for the
reactive power setpoints.

\begin{figure}[t]
\begin{center}
\includegraphics[width=9.0cm]{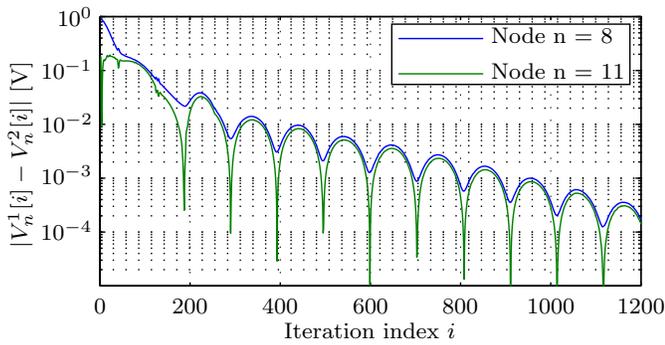}
\end{center}
\vspace{-.3cm}
\caption{Convergence of Algorithm~2: consensus error $|V_{8}^1[i] - V_{8}^2[i]|$ and $|V_{11}^1[i] - V_{11}^2[i]|$ as a function of the ADMM iteration index $i$. }
\label{Fig:OID3results}
\end{figure}

Figure~\ref{Fig:OID3results} represents the discrepancies between local voltages on the line $(8,11)$; specifically, the trajectories of the voltage errors $|V_{8}^1[i] - V_{8}^2[i]|$ and $|V_{11}^1[i] - V_{11}^2[i]|$ are reported  as a function of the ADMM iteration index $i$. The results indicate that the two CEMs consent on the voltage of the branch that connects the two clusters. The ``bumpy'' trend is typical of the ADMM (see e.g.,~\cite{ErsegheADMM, ZhuADMoM }). Similar trajectories were obtained for the inverter setpoints. 

\section{Concluding Remarks} \label{sec:Conclusions}
A suite of decentralized approaches for computing optimal real and
reactive power setpoints for residential photovoltaic (PV) inverters
were developed. The proposed decentralized optimal inverter dispatch
strategy offers a comprehensive framework to share computational
burden and optimization objectives across the distribution network,
while highlighting future business models that will enable customers
to actively participate in distribution-system markets.

\bibliographystyle{IEEEtran}
\bibliography{biblio}

\begin{thebibliography}{10}
\providecommand{\url}[1]{#1}
\csname url@samestyle\endcsname
\providecommand{\newblock}{\relax}
\providecommand{\bibinfo}[2]{#2}
\providecommand{\BIBentrySTDinterwordspacing}{\spaceskip=0pt\relax}
\providecommand{\BIBentryALTinterwordstretchfactor}{4}
\providecommand{\BIBentryALTinterwordspacing}{\spaceskip=\fontdimen2\font plus
\BIBentryALTinterwordstretchfactor\fontdimen3\font minus
  \fontdimen4\font\relax}
\providecommand{\BIBforeignlanguage}[2]{{%
\expandafter\ifx\csname l@#1\endcsname\relax
\typeout{** WARNING: IEEEtran.bst: No hyphenation pattern has been}%
\typeout{** loaded for the language `#1'. Using the pattern for}%
\typeout{** the default language instead.}%
\else
\language=\csname l@#1\endcsname
\fi
#2}}
\providecommand{\BIBdecl}{\relax}
\BIBdecl

\bibitem{Liu08}
E.~Liu and J.~Bebic, ``Distribution system voltage performance analysis for
  high-penetration photovoltaics,'' Feb. 2008, {NREL} Technical Monitor: B.
  Kroposki. Subcontract Report {NREL/SR}-581-42298.

\bibitem{Caramanis_PESTD14}
E.~Ntakou and M.~C. Caramanis, ``Price discovery in dynamic power markets with
  low-voltage distribution-network participants,'' in \emph{IEEE PES Trans. \&
  Distr. Conf.}, Chicago, IL, 2014.

\bibitem{Turitsyn11}
K.~Turitsyn, P.~S\'ulc, S.~Backhaus, and M.~Chertkov, ``Options for control of
  reactive power by distributed photovoltaic generators,'' \emph{Proc. of the
  IEEE}, vol.~99, no.~6, pp. 1063--1073, 2011.

\bibitem{Chertkov-ADMM13}
P.~S\'ulc, S.~Backhaus, and M.~Chertkov, ``Optimal distributed control of
  reactive power via the alternating direction method of multipliers,'' 2013,
  [Online] Available at: \texttt{http://arxiv.org/pdf/1310.5748v1.pdf}.

\bibitem{Aliprantis13}
P.~Jahangiri and D.~C. Aliprantis, ``Distributed {Volt/VAr} control by {PV}
  inverters,'' \emph{IEEE Trans. Power Syst.}, vol.~28, no.~3, pp. 3429--3439,
  Aug. 2013.

\bibitem{Farivar12}
M.~Farivar, R.~Neal, C.~Clarke, and S.~Low, ``Optimal inverter {VAR} control in
  distribution systems with high {PV} penetration,'' in \emph{IEEE PES General
  Meeting}, San Diego, CA, Jul. 2012.

\bibitem{Bolognani13}
S.~Bolognani and S.~Zampieri, ``A distributed control strategy for reactive
  power compensation in smart microgrids,'' \emph{IEEE Trans. on Autom.
  Control}, vol.~58, no.~11, pp. 2818--2833, 2013.

\bibitem{Tonkoski11}
R.~Tonkoski, L.~A.~C. Lopes, and T.~H.~M. El-Fouly, ``Coordinated active power
  curtailment of grid connected {PV} inverters for overvoltage prevention,''
  \emph{IEEE Trans. on Sust. Energy}, vol.~2, no.~2, pp. 139--147, Apr. 2011.

\bibitem{Samadi14}
A.~Samadi, R.~Eriksson, L.~Soder, B.~G. Rawn, and J.~C. Boemer, ``Coordinated
  active power-dependent voltage regulation in distribution grids with pv
  systems,'' \emph{IEEE Trans. on Power Del.}, vol.~29, no.~3, pp. 1454--1464,
  June 2014.

\bibitem{kerstingbook}
W.~H. Kersting, \emph{Distribution System Modeling and Analysis}.\hskip 1em
  plus 0.5em minus 0.4em\relax 2nd ed., Boca Raton, {FL}: {CRC} Press, 2007.

\bibitem{OID}
E.~Dall'Anese, S.~V. Dhople, and G.~B. Giannakis, ``Optimal dispatch of
  photovoltaic inverters in residential distribution systems,'' \emph{IEEE
  Trans. Sustainable Energy}, vol.~5, no.~2, pp. 487--497, Apr. 2014.

\bibitem{Dallanese-TSG13}
E.~Dall'Anese, H.~Zhu, and G.~B. Giannakis, ``Distributed optimal power flow
  for smart microgrids,'' \emph{IEEE Trans. Smart Grid}, vol.~4, no.~3, pp.
  1464--1475, Sep. 2013.

\bibitem{Samadi-SGComm10}
P.~Samadi, A.~H. Mohsenian-Rad, R.~Schober, V.~Wong, and J.~Jatskevich,
  ``Optimal real-time pricing algorithm based on utility maximization for smart
  grid,'' in \emph{Proc. of IEEE Intl. Conf. on Smart Grid Comm.},
  Gaithersburg, MD, 2010.

\bibitem{GatsisTSG12}
N.~Gatsis and G.~B. Giannakis, ``Residential load control: Distributed
  scheduling and convergence with lost {AMI} messages,'' \emph{IEEE Trans.
  Smart Grid}, vol.~3, no.~2, pp. 770--786, 2012.

\bibitem{BeT89}
D.~P. Bertsekas and J.~N. Tsitsiklis, \emph{Parallel and Distributed
  Computation: Numerical Methods}.\hskip 1em plus 0.5em minus 0.4em\relax
  Englewood Cliffs, NJ: Prentice Hall, 1989.

\bibitem{BoydADMoM}
S.~Boyd, N.~Parikh, E.~Chu, B.~Peleato, and J.~Eckstein, ``Distributed
  optimization and statistical learning via the alternating direction method of
  multipliers,'' \emph{Foundations and Trends in Machine Learning}, vol.~3, pp.
  1--122, 2011.

\bibitem{Baldick99}
R.~Baldick, B.~H. Kim, C.~Chase, and Y.~Luo, ``A fast distributed
  implementation of optimal power flow,'' \emph{IEEE Trans. Power Syst.},
  vol.~14, no.~3, pp. 858--864, Aug. 1999.

\bibitem{Nogales03}
F.~J. Nogales, F.~J. Prieto, and A.~J. Conejo, ``A decomposition methodology
  applied to the multi-area optimal power flow problem,'' \emph{Ann. Oper.
  Res.}, no. 120, pp. 99--116, 2003.

\bibitem{Hug09}
G.~Hug-Glanzmann and G.~Andersson, ``Decentralized optimal power flow control
  for overlapping areas in power systems,'' \emph{IEEE Trans. Power Syst.},
  vol.~24, no.~1, pp. 327--336, Feb. 2009.

\bibitem{ErsegheOPFADMM}
T.~Erseghe, ``Distributed optimal power flow using {ADMM},'' \emph{IEEE Trans.
  Power Syst.}, 2014, to appear.

\bibitem{Magnusson14}
S.~Magnusson, P.~C. Weeraddana, and C.~Fischione, ``A distributed approach for
  the optimal power flow problem based on {ADMM} and sequential convex
  approximations,'' Jan. 2014, [Online] Available at:
  \texttt{http://arxiv.org/abs/1401.4621}.

\bibitem{Tse12}
A.~Y. Lam, B.~Zhang, A.~Dom\'{i}nguez-Garc\'{i}a, and D.~Tse, ``Optimal
  distributed voltage regulation in power distribution networks,'' 2012,
  [Online] Available at \texttt{http://arxiv.org/abs/1204.5226v1}.

\bibitem{Zhu_DSE}
H.~Zhu and G.~B. Giannakis, ``Multi-area state estimation using distributed
  {SDP} for nonlinear power systems,'' in \emph{3rd Int. Conf. Smart Grid
  Comm.}, Tainan City, Taiwan, Nov. 2012.

\bibitem{Kekatos_stateest13}
V.~Kekatos and G.~B. Giannakis, ``Distributed robust power system state
  estimation,'' \emph{IEEE Trans. Power Syst.}, vol.~28, no.~2, pp. 1617--1626,
  May 2013.

\bibitem{DallAneseTSTE14}
E.~Dall'Anese and G.~B. Giannakis, ``Risk-constrained microgrid reconfiguration
  using group sparsity,'' May 2014, to appear; see also:
  \texttt{http://arxiv.org/abs/1306.1820}.

\bibitem{Wiesel11}
A.~T. Puig, A.~Wiesel, G.~Fleury, and A.~O. Hero, ``Multidimensional
  shrinkage-thresholding operator and group {LASSO} penalties,'' \emph{IEEE
  Sig. Proc. Letters}, vol.~18, no.~6, pp. 363--366, Jun. 2011.

\bibitem{LavaeiLow}
J.~Lavaei and S.~H. Low, ``Zero duality gap in optimal power flow problem,''
  \emph{IEEE Trans. Power Syst.}, vol.~1, no.~1, pp. 92--107, 2012.

\bibitem{Lavaei_tree}
J.~Lavaei, D.~Tse, and B.~Zhang, ``Geometry of power flows and optimization in
  distribution networks,'' in \emph{IEEE PES General Meeting}, San Diego, CA,
  2012.

\bibitem{ErsegheADMM}
T.~Erseghe, D.~Zennaro, E.~Dall'Anese, and L.~Vangelista, ``Fast consensus by
  the alternating direction multipliers method,'' \emph{IEEE Trans. Sig.
  Proc.}, vol.~59, no.~11, pp. 5523--5537, 2011.

\bibitem{Vandenberghe96}
L.~Vandenberghe and S.~Boyd, ``Semidefinite programming,'' \emph{{SIAM}
  Review}, vol.~38, no.~1, pp. 49--95, Mar. 1996.

\bibitem{Grone84}
R.~Grone, C.~R. Johnson, E.~M. S\'a, and H.~Wolkowicz, ``Positive definite
  completions of partial {H}ermitian matrices,'' \emph{Linear Algebra and its
  Applications}, vol.~58, pp. 109--124, 1984.

\bibitem{Jabr12}
R.~A. Jabr, ``Exploiting sparsity in {SDP} relaxations of the {OPF} problem,''
  \emph{IEEE Trans. Power Syst.}, vol.~2, no.~27, pp. 1138--1139, May 2012.

\bibitem{Braun10}
M.~Braun, J.~K\"unschner, T.~Stetz, and B.~Engel, ``Cost optimal sizing of
  photovoltaic inverters -- influence of of new grid codes and cost
  reductions,'' in \emph{Proc. of 25th Europ. {PV} Solar Energy Conf. and
  Exhib.}, Valencia, Spain, Sep. 2010.

\bibitem{ZhuADMoM}
H.~Zhu, G.~B. Giannakis, and A.~Cano, ``Distributed in-network channel
  decoding,'' \emph{IEEE Trans. Sig. Proc.}, vol.~57, no.~10, pp. 3970--3983,
  2009.

\end{thebibliography}

\end{document}